\documentclass[
    onefignum,
    onetabnum,
]{siamart171218}

\usepackage[%
left=1.5in,
right=1.5in,
bottom=0.75in,
top=1in,
footskip=0.5in,
]{geometry}

\usepackage{graphicx}
\usepackage{epstopdf}
\usepackage{algorithm} %provides pseudocode framework
\usepackage{algorithmic} 
\usepackage[font=footnotesize]{caption}
\usepackage{subcaption}
\usepackage{amsfonts,amsmath,amssymb,mathtools} 
\usepackage[mathscr]{eucal}
\usepackage{graphicx}
\usepackage{bbm} %for indicator function
\usepackage{enumitem}
\usepackage{cite}
\usepackage{placeins}
\usepackage{csquotes}
\ifpdf
\DeclareGraphicsExtensions{.eps,.pdf,.png,.jpg}
\else
\DeclareGraphicsExtensions{.eps}
\fi

\title{Controlling a Vlasov--Poisson plasma \\ by a Particle-In-Cell method \\ based on a Monte Carlo framework}

\headers{Controlling plasma in a Monte Carlo framework}{J. Bartsch, P. Knopf, S. Scheurer, J. Weber}

\author{Jan Bartsch\thanks{Department of Mathematics, University of Konstanz, Germany (\email{jan.bartsch@uni-konstanz.de}).}
	\and Patrik Knopf\thanks{Department of Mathematics, University of Regensburg, Germany
		(\email{patrik.knopf@ur.de}).}
	\and Stefania Scheurer\thanks{SimTech Cluster, University of Stuttgart, Germany
		(\email{stefania.scheurer@iws.uni-stuttgart.de}).}
	\and  Jörg Weber\thanks{Faculty of Mathematics, University of Vienna, Austria
		(\email{joerg.weber@univie.ac.at}).}
}

\newtheorem{rem}[theorem]{Remark}

\newtheorem{assumption}[theorem]{Assumption}

\makeatletter
\renewenvironment{algorithm}
{%
    \bigskip
    \begin{center}
    \refstepcounter{algorithm}% New algorithm
        \hrule height.8pt depth0pt \kern2pt% \@fs@pre for \@fs@ruled
        \renewcommand{\caption}[2][\relax]{% Make a new \caption
        {\raggedright\textbf{\fname@algorithm~\thealgorithm} ##2\par}%
        \ifx\relax##1\relax % #1 is \relax
         \addcontentsline{loa}{algorithm}{\protect\numberline{\thealgorithm}##2}%
        \else % #1 is not \relax
         \addcontentsline{loa}{algorithm}{\protect\numberline{\thealgorithm}##1}%
        \fi
        \kern2pt\hrule\kern2pt
        }
}
{%
     \kern2pt\hrule\relax% \@fs@post for \@fs@ruled
     \bigskip
   \end{center}
}
\makeatother

\def\EE{\mathbb{E}}

\def\NN{\mathbb{N}}

\def\RR{\mathbb{R}}
\def\XX{\mathbb{X}}
\def\VV{\mathbb{V}}

\newcommand{\Bad}{\mathbb{B}_{\mathrm{ad}}}

\newcommand{\BigO}[1]{\mathcal{O}(#1)}

\DeclareMathOperator\supp{supp}

\begin{document}
	
\maketitle

\begin{abstract}
The Vlasov--Poisson system describes the time evolution of a plasma in the so-called collisionless regime. The investigation of a high-temperature plasma that is influenced by an exterior magnetic field is one of the most significant aspects of thermonuclear fusion research. In this paper, we formulate and analyze a kinetic optimal control problem for the Vlasov--Poisson system where the control is represented by an external magnetic field.
The main goal of such optimal control problems is to confine the plasma in a certain region in phase space.
We first investigate the optimal control problem in terms of mathematical analysis, i.e., we show the existence of at least one global minimizer and we rigorously derive a first-order necessary optimality condition for local minimizers by the adjoint approach.
Then, we build a Monte Carlo framework to solve the state equations as well as the adjoint equations by means of a Particle-In-Cell method, and we apply a nonlinear conjugate gradient method to solve the optimization problem.
Eventually, we present numerical experiments that successfully validate our optimization framework.
\end{abstract}

% REQUIRED
\begin{keywords}
	Vlasov--Poisson system, ensemble optimal control problems, multi-species plasma, Monte Carlo methods, Particle-In-Cell method, magnetic confinement.
\end{keywords}

% REQUIRED
\begin{AMS}
	49J20, %Existence theories for optimal control problems involving partial differential equations
    35Q83, %Vlasov equations 
    49M05, %Methods based on necessary conditions
    62P35, %Statistics: Applications to physics
    65C05, %Numerical analysis/Simulation: Monte Carlo methods
    90C15  %Stochastic programming
\end{AMS}

\setlength\parskip{0.5ex}

\section{Introduction} 
The investigation of a high-temperature plasma that is influenced by an exterior magnetic field is one of the most important aspects in the research on thermonuclear fusion \cite{Chen1984introductionPlasmaFusion, Stacey2005FusionPlasmaPhysics}.
In this context, the magnetic confinement of plasmas is of particular interest. 
The idea of magnetic confinement is to apply magnetic fields to influence the plasma in such a way that it keeps a certain distance to the reactor wall. This is necessary because a plasma colliding with the wall will cool down (which interrupts the fusion process) and might even damage the reactor.
In this paper, we investigate optimal control problems where a plasma is supposed to be influenced in a desired way. In particular, with regard to magnetic confinement, we want to control a plasma in such a way that its particles do not leave a prescribed region in phase space.

To mathematically describe the time evolution of a plasma consisting of electrons and ions, we use the kinetic theory of gases. 
We assume that there is only one type of ions with mass $m_+$ and charge $q$, and we denote by $m_-$ the mass of the electrons having the charge $-q$.
We consider the spatial domain to be $\RR^{d_x}$ and the velocity space to be $\RR^{d_v}$.
Moreover, to represent the density in phase space, $f^+ = f^+(t,x,v) \geq 0$ denotes the distribution function of the ions, whereas $f^- = f^-(t,x,v) \geq 0$ stands for the distribution function of the electrons. Furthermore, $\rho_f$ stands for the corresponding charge density that is defined as
\begin{align}
\label{def:rho}
	\rho_f(t,x) = \int_{\RR^{d_v}} f^+(t,x,v) \,\mathrm dv 
    - \int_{\RR^{d_v}} f^-(t,x,v)\,\mathrm dv.
\end{align}
In this paper, we consider plasmas in the so-called \textit{collisionless regime}. That is, we expect that all collisions between the particles are so rare that they can be neglected. This is a widely used and legitimate approximation especially for fusion plasmas in a tokamak device (see, e.g., \cite{BadsiHerda2016ModelSimPlasma, Miyamoto2016PlasmaControlledFusion}).
Furthermore, we only take electro-static effects into account. This means that merely the self-induced electric field is considered whereas the self-induced magnetic field is neglected. 
This is justified as in order to control the plasma, we will also apply an external magnetic field $B = B(t,x)$, which is typically much stronger than the self-induced magnetic field. 
In this case, the evolution of the distribution functions $f^\pm$ over time is described by the Vlasov--Poisson system. In this model, the time evolution of the distribution functions $f^\pm$ is described by the following partial differential equations:
\begin{subequations}
	\begin{align}
		\partial_tf^- + v\cdot\partial_xf^- -\frac{q}{m_-}(E_f+v\times B)\cdot\partial_vf^-=0,
		\label{eq:Poisson_equation_electrons_unscaled}\\
		\partial_tf^+ + v\cdot\partial_xf^+ +\frac{q}{m_+}(E_f+ v\times B)\cdot\partial_vf^+=0.
		\label{eq:Poisson_equation_ions_unscaled}
	\end{align}
\end{subequations}
These first-order PDEs are referred to as the \textit{Vlasov equations},
and $E_f$ is the self-consistent electric field induced by the plasma particles. As it is a common notation in the context of Vlasov equations, we write $\partial_x \coloneqq  \nabla_x$ and $\partial_v \coloneqq  \nabla_v$ to denote the partial gradients with respect to $x$ and $v$, respectively.

Due to the presence of electrons and ions, the dynamics of \eqref{eq:Poisson_equation_electrons_unscaled} and \eqref{eq:Poisson_equation_ions_unscaled} are very heterogeneous.
In particular, we have several spatio-temporal scales of different orders of magnitudes that lead to the fact that the numerical simulation of plasmas is very challenging. 
In order to keep our discussion physically relevant and to ensure that the system parameters are chosen in reasonable orders of magnitudes, we employ a scaling by the characteristic quantities associated with our plasma.
These are the \textit{plasma frequency} $\omega_{\text{pe}}$, 
the \textit{thermal velocities} $v_{\text{th},\pm}$ and the \textit{Debye length} $L_{\text{DE}}$. They can be expressed by the following formulas:
\begin{align*}
	 \omega_{\text{pe}} = \sqrt{\frac{n_-q^2}{\varepsilon_0m_-}}, \qquad\qquad v_{\text{th},\pm} = \sqrt{\frac{k_BT_\pm}{m_\pm}}, \qquad\qquad
	 L_{\text{DE}} = \frac{v_{\text{th},-}}{\omega_{\text{pe}}}.
\end{align*} 
Here, $n_-$ is the \textit{number density} of electrons, $\varepsilon_0$ is the \textit{dielectric constant}, $T_{\pm}$ are the \textit{characteristic temperatures} of the corresponding species and $k_B$ is the \textit{Boltzmann constant}.
Using these definitions, we express the time in units of $\omega_{\text{pe}}^{-1}$, 
the position in units of $L_{\text{DE}}$ and the velocity in units of $v_{\text{th},\pm}$.
To account for the different scaling of the velocities, we introduce the constant dimensionless factors
\begin{align*}
	\mu_x^+=\frac{v_{\text{th},+}}{v_{\text{th},-}},
	\qquad
	\mu_v^+=\frac{m_-}{\mu_x^+m_+},
    \qquad
    \mu_x^- = 1,
    \qquad \mu_v^- = -1.
\end{align*}
These scaling parameters are very useful as they make it possible to study the time evolutions of the electrons and of the ions on the same time scale.
By means of this rescaling, the Vlasov--Poisson equations \eqref{eq:Poisson_equation_electrons_unscaled} and \eqref{eq:Poisson_equation_ions_unscaled} (for prescribed initial data $\mathring f = (\mathring f^-,\mathring f^+)$) can be reformulated as 
\begin{subequations}
    \label{VP:PM}
	\begin{align}
    	&\partial_t f^\pm + \mu_x^\pm v \cdot\partial_xf^\pm + \mu_v^\pm (E_f+ \mu_x^\pm v\times B) \cdot \partial_vf^\pm=0,
        \\
		&f^\pm\vert_{t=0} = \mathring f^\pm
	\end{align}
\end{subequations}
(cf.~\cite[Appx.~I]{MangETAL02}).
Here, depending on the dimension $d_x$, the self-consistent electric field is given by
\begin{align}\label{eq:defE}
        E_f(t,x) = 
        \begin{cases}
            \int_{-\infty}^x \rho_f(t,y) \, \mathrm dy 
            &\text{if $d_x=1$}, \\[1ex]
            \frac{1}{2\pi} \int_{\mathbb R^2} \frac{(x-y)}{|x-y|^2} \rho_f(t,y) \, \mathrm dy 
            &\text{if $d_x=2$}, \\[1ex]
            \frac{1}{4\pi} \int_{\mathbb R^3} \frac{(x-y)}{|x-y|^3} \rho_f(t,y) \, \mathrm dy 
            &\text{if $d_x=3$},
        \end{cases}
\end{align}
where $\rho_f$ is defined as in \eqref{def:rho}. 

Typically, the parameters $T_\pm$ and $n_-$ are of order $10^7$--$10^8 \,\text{K}$ and $10^{20}$--$10^{21}\,\text{m}^{-3}$, respectively, so that $\mu_x^+$ and $\mu_v^+$ are both of order $10^{-2}$. To get a rough flavor of the appearing scales, time, length, electron velocity, ion velocity, electric field, and magnetic field are measured in units of
\begin{gather*}
    \omega_{\text{pe}}^{-1}\sim 10^{-12}\text{--}10^{-11}\,\text{s},\qquad L_{\text{DE}}\sim 10^{-4}\,\text{m},\qquad v_{\text{th},-}\sim 10^7\text{--}10^8\,\frac{\text{m}}{\text{s}},\qquad v_{\text{th},+}\sim 10^6\,\frac{\text{m}}{\text{s}},\\
    E^0=\frac{m_-L_{\text{DE}}\omega_{\text{pe}}^2}{q}\sim 10^8\text{--}10^9\,\frac{\text{V}}{\text{m}},\qquad 
    B^0=\frac{m_-\omega_{\text{pe}}}{q}\sim 10^0\text{--}10^1\,\text{T}.
\end{gather*}
(Recall that $1\,\text{V} = 1\,\text{kg}\cdot\text{m}^2\cdot(\text{C}\cdot\text{s}^2)^{-1}$ and $1\,\text{T} = 1\,\text{kg}\cdot(\text{C}\cdot\text{s})^{-1}$.)
Typical external magnetic fields applied in real-world fusion reactors are roughly of the size of $B^0$.
In fact, it is technically possible to generate magnetic fields that exceed $45\,\text{T}$, cf.~\cite{Hahn2019TeslaMagneticField}. Therefore, the modulus of our dedimensionalized magnetic field $B$ should roughly be of size $1$.  
Moreover, the dedimensionalized distribution functions $f^\pm$ are measured in the unit $n_-v_{\text{th},\pm}^{-3}$.

Furthermore, we assume that the initial total charge vanishes, i.e.,\[\int_{\RR^{d_x}}\rho_{\mathring f}(x)\, \mathrm dx=
\int_{\RR^{d_x}}\int_{\RR^{d_v}}(\mathring f^+-\mathring f^-)(x,v)\, \mathrm dv\mathrm dx=0,\]
From a physical point of view, this is a reasonable assumption as 
non-ionized atoms possess the same number of protons as electrons. We point out that the Vlasov equations preserve this charge neutrality over the course of time.
Notice that there is no need to impose boundary conditions for $f^\pm$, since we do not prescribe any bounds for the particles (at least in the theoretical part of the paper) and the distribution functions $f^\pm$ will be compactly supported in phase space for all times.

Even though the Vlasov--Poisson system is usually considered in the full-dimensional setting $d_x=d_v=3$, we will also consider its lower dimensional versions.
They describe the situation in which the distribution functions $f^\pm$ are constant with respect to certain space or velocity components. In the following, we assume $d_v \in \{2,3\}$ for the velocity dimension and $d_x\in \mathbb N$ with $d_v-1\le d_x \le d_v$ for the spatial dimension. We write $d=(d_x,d_v)$. Then the external magnetic field is a function $B\colon[0,\infty)\times \RR^{d_x} \to \RR^{d_B}$ with
\begin{align} \label{def:db}
    d_B \coloneqq  
    \begin{cases}
        1 &\text{if $d_v=2$},\\
        3 &\text{if $d_v=3$}.
    \end{cases}
\end{align}
This is because if the velocities are merely two-dimensional, we can only consider magnetic fields that are orthogonal to the $(v_1,v_2)$-plane and thus, only one vector component of $B$ needs to be specified.
Moreover, the internal electric field is a function $E_f\colon[0,\infty)\times \RR^{d_x} \to \RR^{d_x}$. To formulate the $d$-dimensional version of \eqref{VP:PM}, we introduce the following notation:
\begin{align*}
    F \otimes G \coloneqq  
    \begin{cases}
        F \times G &\text{if $F,G\in\RR^3$},\\
        \begin{psmallmatrix} F_2 \\ - F_1 \end{psmallmatrix} G &\text{if $F\in\RR^2$, $G\in\RR$}.
    \end{cases}
\end{align*}
Moreover, for any vectors $F\in \RR^{n}$, $G\in\RR^m$ with $n,m\in\mathbb N$, we interpret their sum and their scalar product in such a way that the shorter vector is filled up by zeros. This means
\begin{align*}
    F + G \coloneqq  
    \begin{cases}
        F + G &\text{if $n=m$},\\
        \begin{psmallmatrix} F \\ 0 \end{psmallmatrix} + G &\text{if $n<m$},\\
        F + \begin{psmallmatrix} G \\ 0 \end{psmallmatrix} &\text{if $n>m$},
    \end{cases}
    \qquad
    F \cdot G \coloneqq  
    \begin{cases}
        F \cdot G &\text{if $n=m$},\\
        \begin{psmallmatrix} F \\ 0 \end{psmallmatrix} \cdot G &\text{if $n<m$},\\
        F\cdot \begin{psmallmatrix} G \\ 0 \end{psmallmatrix} &\text{if $n>m$}.
    \end{cases}
\end{align*}
In this way, together with \eqref{eq:defE}, we can formulate the $d$-dimensional version of the Vlasov--Poisson system \eqref{VP:PM} as
\begin{subequations}
\label{VP:PM:d}
\begin{align}
    &\partial_t f^\pm + \mu_x^\pm v \cdot\partial_x f^\pm + \mu_v^\pm (E_f + \mu_x^\pm v\otimes B) \cdot \partial_v f^\pm=0,
    \\
    &f^\pm\vert_{t=0} = \mathring f^\pm.
\end{align}
\end{subequations}

The Vlasov--Poisson system has already been investigated in the literature from many different viewpoints. 
A first local classical well-posedness result for the standard Vlasov--Poisson system (i.e., with $B\equiv 0$) for compactly supported, continuously differential initial data was established by Kurth~\cite{Kurth1952}. Later a continuation criterion stating that such local-in-time solutions can be extended as long as their velocity/momentum support is under control was proven by Batt~\cite{Batt1977}. 
In the case $d=(2,2)$ and $B\equiv 0$, the existence of a unique global classical solution was shown by Ukai \& Okabe~\cite{Ukai1978} but also (almost simultaneously) by Wollman~\cite{Wollman1980}. 
In the full-dimensional case $d=(3,3)$ (still with $B\equiv 0$), the global existence of a unique global classical solution was established independently and simultaneously in the seminal works by Pfaffelmoser \cite{Pfaffelmoser1992} and by Lions \& Perthame \cite{Lions1991}. It is worth mentioning that both proofs differ greatly with respect to the employed techniques. A greatly simplified version of Pfaffelmoser's proof was given by Schaeffer in \cite{Schaeffer1991}. Later, the Pfaffelmoser/Schaeffer proof was adapted by Knopf~\cite{Knopf_Diss,Knopf2018OCP_VP} to prove global existence in dimension $d=(3,3)$ also for a class of non-trivial external magnetic fields. An analogous global existence result for non-trivial magnetic fields in the case $d=(2,3)$ was established by Knopf \& Weber~\cite{KnopfWeber2021VlasovPoissonSteadyStates}. Moreover, for the investigation of steady states (i.e., time-independent solutions) of the Vlasov--Poisson system under the influence of an external magnetic field, we refer to \cite{Knopf2019ConfinedStatesCylinder,KnopfWeber2021VlasovPoissonSteadyStates,Rein1992,Braasch1999,Belyaeva2021,Skubachevskii2014} as well as the references therein. Moreover, the problem of magnetic confinement, also for slightly different plasma models, was studied, e.g., in \cite{CCM12,Han10,NNS15,Zhang21}.

Furthermore, optimal control problems for the Vlasov--Poisson system, where an external magnetic field represents the control, were investigated by Knopf~\cite{Knopf_Diss,Knopf2018OCP_VP} as well as by Knopf \& Weber~\cite{KnopfWeber2020_OCPVP_Coils}. Afterwards, optimal control problems for the more complicated Vlasov--Maxwell equations were studied by Weber~\cite{W20:OCRVM,W21:OCRVM2D}. We further want to mention that the controllability of the Vlasov--Poisson system (in the context of control theory) was investigated by Glass~\cite{Glass2003} as well as Glass \& Han-Kwan~\cite{Glass2012}. 

It the present work, we extend the optimal control theory developed in \cite{Knopf_Diss,Knopf2018OCP_VP} to make it more suitable for its numerical investigation. In particular, we show that the required regularity of admissible controls can be relaxed. This leads to a simpler first-order optimality condition that can be handled more easily by numerical methods.
To solve the optimality system, we use a Particle-In-Cell method based on a Monte Carlo framework developed by Bartsch~\cite{Bartsch2021MOCOKI, Bartsch2020OCPKS}.

This paper is organized as follows: In Section~\ref{sec:prelim}, we first introduce some notation and preliminaries. Afterwards, in Subsection~\ref{sec:WellPosedness_classical} and Subsection~\ref{sec:WellPosedness_strong}, we discuss the well-posedness of the Vlasov--Poisson system for magnetic fields of class $C^1$ and $W^{1,\infty}$, respectively. This allows for the definition of a control-to-state operator whose most important properties are shown in Subsection~\ref{sec:ControlToStateOperator}. In Subsection~\ref{sec:AdmissibleControls}, we introduce our set of admissible controls, and in Subsection~\ref{sec:Formulation_OPC}, we formulate our kinetic optimal control-in-the-force problem in the framework of ensemble control problems. In particular, we illustrate the construction of our cost functional including ensemble trajectory terms and ﬁnal-conﬁguration contributions as well as a penalization term for the control. 
For this setting, we then prove the existence of solutions in the given control space.
Afterwards, in Section \ref{sec:OptimalitySystem}, we use the adjoint approach to formulate the first-order necessary optimality condition for locally optimal controls. In Section \ref{sec:MC_framework}, we present a detailed discussion of the numerical methods that we will use to solve the optimality system. 
In particular, we explain how to compute the gradient of the cost-functional, which is then used to implement a nonlinear conjugate gradient scheme (see, e.g., \cite{BorziSchulz2011}).
In Section \ref{sec:NumericalExperiments}, we present the results of numerical experiments that successfully validate our implementation. 
For our numerical solver of the Vlasov--Poisson system, we use the famous \textit{Landau damping} effect as well as the \textit{two-stream instability} effect as test cases. 
We further illustrate the functionality of our optimization method by solving an optimal control problem for magnetic confinement of a plasma.
In contrast to the uncontrolled case, we actually observe confinement of the controlled plasma during the prescribed time interval.

\section{Notation and preliminaries}\label{sec:prelim}

Throughout, we shall write $z=(x,v)\in\RR^{d_x}\times\RR^{d_v}$ for brevity whenever convenient.

For any tuple $f=(f^-,f^+)$ of two functions, norms of $f$ are always understood to be sum of the norms of the two components, that is,
\[\|f\|_X\coloneqq\|f^-\|_X+\|f^+\|_X\]
for any norm $\|\cdot\|_X$.

For any function $\varrho\in L^1(\mathbb{R}^{d_x})$ with compact support and zero mean, we shall write $(-\Delta_x)^{-1} \varrho = U_\varrho$ to denote the solution of Poisson's equation
\begin{gather*}
    -U_\varrho'' = \varrho
    \quad\text{with}\quad
    \lim_{x\to-\infty}U_\varrho(x) = 0  
    \quad\text{and}\quad
    \lim_{x\to\infty}U_\varrho(x) \text{ exists in }\RR 
\end{gather*}
if $d_x = 1$ or
\begin{gather*}
    -\Delta_x U_\varrho = \varrho
    \quad\text{with}\quad
    U_\varrho(x) \to 0 \quad\text{as }|x|\to \infty
\end{gather*}
if $d_x\in\{2,3\}$, respectively.

We first present some important estimates, which will be used later to bound the self-consistent electric field in terms of the underlying charge density.

\begin{lemma} \label{lem:E}
    Let $d_x\in \{1,2,3\}$ and let $\varrho \in L^2(\mathbb{R}^{d_x})$ with $\mathrm{supp}\, \varrho \subset B_r(\mathbb{R}^{d_x})$ for some $r>0$ and $\int_{\RR^{d_x}}\varrho(x)\, \mathrm dx=0$. Moreover, we define 
    \begin{align*}
        E(x) \coloneqq  -\nabla_x (-\Delta_x)^{-1} \varrho (x) =
        \begin{cases}
            \int_{-\infty}^x \varrho(y) \, \mathrm dy 
            &\text{if $d_x=1$}, \\[1ex]
            \frac{1}{2\pi} \int_{\mathbb R^2} \frac{(x-y)}{|x-y|^2} \varrho(y) \, \mathrm dy 
            &\text{if $d_x=2$}, \\[1ex]
            \frac{1}{4\pi} \int_{\mathbb R^3} \frac{(x-y)}{|x-y|^3} \varrho(y) \, \mathrm dy 
            &\text{if $d_x=3$}.
        \end{cases}
    \end{align*}
    for all $x\in\mathbb R^{d_x}$.
    Then the following holds:
    \begin{enumerate}[label = \textnormal{(\alph*)}, leftmargin=*]
        \item\label{est:E:L2} $E \in L^2(\mathbb R^{d_x})$, and there exists a constant $C = C(d_x,r) >0$ such that
        \begin{align*}
            \|E\|_{L^2(\RR^{d_x})} \le C \|\varrho\|_{L^2(\RR^{d_x})}.
        \end{align*}
        \item\label{est:E:Linf} 
        If $d_x=1$, then $E \in L^\infty(\mathbb R^{d_x})$ with
        \begin{align*}
            \|E\|_{L^\infty(\RR^{d_x})} \le \|\varrho\|_{L^1(\RR^{d_x})}.
        \end{align*}
        If $d_x\in\{2,3\}$ and additionally $\varrho \in L^\infty(\mathbb R^{d_x})$, then $E \in L^\infty(\mathbb R^{d_x})$ with
        \begin{align*}
            \|E\|_{L^\infty(\RR^{d_x})} \le C \|\varrho\|_\infty
        \end{align*}
        for some constant $C = C(d_x,r) >0$.
        \item\label{est:E:W1inf} If additionally $\varrho \in W^{1,\infty}(\mathbb R^{d_x})$, then $E \in W^{1,\infty}(\mathbb R^{d_x})$ with
        \begin{align*}
            \|E\|_{W^{1,\infty}(\RR^{d_x})} 
            \le c \Big[ \big(1 + \|\varrho\|_\infty\big) \big( 1 + \ln_+\|\nabla \varrho \|_\infty\big) + \|\varrho\|_1 \Big]
        \end{align*}
        for some constant $c = c(d_x,r) >0$.
    \end{enumerate}
\end{lemma}

\begin{proof}
    In the case $d_x = 1$, the proof is straightforward.
    For $d_x = 3$, \ref{est:E:L2} and \ref{est:E:Linf} were presented, e.g., in \cite[Lemma~2(c)]{Knopf2018OCP_VP}, and \ref{est:E:W1inf} was shown in \cite[Lemma~P1]{Rein2007Review}. The case $d_x=2$ can be established analogously to the case $d_x = 3$.
\end{proof}

We point out that at several points in the above lemma, the assumptions on $\varrho$ could be weakened and the estimates on $E$ could be improved in the case $d_x\in\{1,2\}$. However, the results are sufficient for later usage and allow us to handle all cases $d_x\in\{1,2,3\}$ simultaneously.

\section{Analysis of the optimal control problem}\label{sec:Analysis_OCP}
In this section, we consider a finite time horizon $[0,T]$ with $T>0$ and the phase space $\RR^{D}$ with $D\coloneqq  d_x+d_v$, where 
$d\coloneqq (d_x,d_v) \in \{(1,2),(2,2),(2,3)\}$, 
on which we investigate the evolution of the distribution functions $f^\pm$ governed by \eqref{VP:PM:d}. 
For arbitrary, fixed radii $R_x, R_v>0$, we define 
\[
    \Omega_x \coloneqq B_{R_x}^{d_x}(0), \quad \Omega_v \coloneqq B_{R_v}^{d_v}(0), \quad \Omega\coloneqq  \Omega_x\times\Omega_v,
\]
where the superscript indicates the dimension of the ball.
We further fix arbitrary initial data $\mathring{f} = (\mathring{f}^-, \mathring{f}^+)$ with $\mathring{f}^\pm \in C^1_c(\RR^D;\RR^+_0)$ and $\supp \mathring{f}^\pm \subset \Omega$. 
We will see that the final time $T>0$ can be chosen in such a way that the support of the solution remains contained in $\Omega$. In this context, it is crucial that this choice of $T$ does not depend on the choice of the external magnetic field. 
As a consequence, it suffices to define the external magnetic field merely on $[0,T]\times \overline{\Omega_x}$. Hence, it is a function $B:[0,T]\times\overline{\Omega_x} \to \RR^{d_B}$ where $d_B$ is given by \eqref{def:db}.

\subsection{Classical well-posedness for magnetic fields of class \texorpdfstring{$C^1$}{C\^{}1}}
\label{sec:WellPosedness_classical}

We first prove classical well-posedness of the Vlasov--Poisson system \eqref{VP:PM:d} on $[0,T]$ for continuously differentiable external magnetic fields. In particular, the final time $T$ can be chosen independent of $B$ such that the support of $f^\pm$ remains contained in $\overline\Omega$.

\begin{proposition}[Classical well-posedness for magnetic fields of class $C^1$]
\label{prop:existence_classical_solutions}
    There exists a final time $T>0$ depending only on $\mathring{f}$, $R_x$, $R_v$ and $d$ such that 
	for all magnetic fields $B \in C^1([0,T] \times \overline{\Omega_x};\RR^{d_B})$, 
	  there exists a unique classical solution $f=(f^-,f^+)$ with $f^\pm \in C^1([0,T]\times\RR^D)$ of the Vlasov--Poisson system \eqref{VP:PM:d} to the
    initial data $\mathring f$ with
	$\supp f^\pm(t) \subset \overline\Omega$ for all $t \in [0,T]$.
\end{proposition}

\begin{proof}
    First, let $T>0$ be arbitrary.
    Since $B \in C^1([0,T] \times \overline{\Omega_x};\RR^{d_B})$,
    it can be extended to a function $\overline{B} \in C^1_b\big(\RR\times\RR^{d_x};\RR^{d_B})$ meaning that 
    $\overline{B}\vert_{[0,T]\times\overline{\Omega_x}} = B$. 

    Let us first assume that $f$ is a classical solution of the Vlasov--Poisson system \eqref{VP:PM:d}.
    We write $Z^\pm=(X^\pm,V^\pm)(s,t,z)$ to denote the unique solutions of the characteristic systems
    \begin{align}
        \label{charsys}
        \dot x(s) = v(s),\quad \dot v(s) =  \mu_v^\pm E_f\big(s,x(s)\big) + \mu_x^\pm \mu_v^\pm v(s) \otimes \overline{B}\big(s,x(s)\big)
    \end{align}
    satisfying the initial conditions $Z^\pm(t,t,z) = z$, respectively, for any given $z=(x,v)\in\mathbb R^6$.  Since $\supp \mathring{f} \subset \Omega = \Omega_x\times\Omega_v$, there exist $R_x^0 \,\in\, ]0,R_x[$ and $R_v^0 \,\in\, ]0,R_v[$ depending only on $\mathring{f}$ such that $\supp \mathring{f} \subset B_{R_x^0}(0) \times B_{R_v^0}(0)$.
    To bound the supports of $f^\pm$, we define the quantities%
    \begin{subequations}
    \begin{align}
        \label{def:X}
        \XX(t) &\coloneqq  \underset{\pm}{\max} \Big\{ |x| \;\Big\vert\; \exists v\in\RR^{d_v},\, s\in[0,t]: (x,v) \in \supp f^\pm(s) \Big\},\\
        \label{def:V}
        \VV(t) &\coloneqq  \underset{\pm}{\max} \Big\{ |v| \;\Big\vert\; \exists x\in\RR^{d_x},\, s\in[0,t]: (x,v) \in \supp f^\pm(s) \Big\}
    \end{align}
    \end{subequations}
    for all $t$ as long as the solution exists. As $f$ is constant along the characteristic flow, these quantities can alternatively be expressed as%
    \begin{align}
        \label{def:X:char}
        \XX(t) 
        &\coloneqq  \underset{\pm}{\max} \Big\{ |X(s,0,x,v)| \;\Big\vert\; s\in[0,t] \wedge \exists v\in\RR^{d_v}: (x,v) \in \supp \mathring{f}\Big\}, \\
        \label{def:V:char}
        \VV(t) 
        &\coloneqq  \underset{\pm}{\max} \Big\{ |V(s,0,x,v)| \;\Big\vert\; s\in[0,t] \wedge \exists x\in\RR^{d_x}: (x,v) \in \supp \mathring{f}\Big\}.
    \end{align}
    In the following, we write $C$ to denote generic positive constants depending only on $\mathring{f}$ and $d$ which are allowed to change their value from line to line.
    
    \textit{The case $d=(2,3)$.}
    The existence of a unique global classical solution to the magnetic field $\overline B$ was established in \cite[Theorem~3.1]{KnopfWeber2021VlasovPoissonSteadyStates}. Moreover, it was shown in the proof of this theorem that 
    \begin{align}
        \label{EST:X23}
        \XX(t) &\le R_x^0 + \int_0^t \VV(s) \,\mathrm ds \le R_x^0 + t\, \VV(t),\\
        \label{EST:V23}
        \VV(t) &\le \left[ (1+R_v^0)^{\frac 14} + C \int_0^t (1+s)^{\frac 14} \,\mathrm ds \right]^4 - 1 
    \end{align}  
    for all $t\ge 0$. Note that the right-hand sides of \eqref{EST:X23} and \eqref{EST:V23} are increasing and attain the values $R_x^0$ and $R_v^0$, respectively, at $t=0$. Thus, by choosing $T$ sufficiently small, we ensure $\XX(t) < R_x$ and $\VV(t) < R_v$
    for all $t\in [0,T]$. We point out that this choice of $T$ depends only on $\mathring{f}$ and $d$ since $R_x^0$, $R_v^0$ and $C$ depend at most on $\mathring{f}$ and $d$. Consequently, the solution $f$ depends only on $\overline{B}\vert_{[0,T]\times\overline{\Omega_x}} = B$ but not on the extension $\overline{B}$.

    \textit{The case $d=(2,2)$.} Here, one can follow the same ideas as in the case $d=(2,3)$. In fact, things are even easier, so we omit the details.

    \textit{The case $d=(1,2)$.} This case is even easier than the higher dimensional cases, but since we later study this case numerically, we provide the details. As in the higher dimensional cases, the existence of a local classical solution $f$ to the magnetic field $\overline B$ existing on its right-maximal time interval $[0,T_{\max}[$ can be established by a standard fixed-point argument (see, e.g., \cite[Section~1.2]{Rein2007Review}). 
    Moreover, it is well-known that in order to prove $T_{\max}=\infty$, it suffices to show that $\VV(t)$ remains bounded on any finite subinterval of $[0,T_{\max}[$ (see, e.g., Theorem~\cite[Theorem~2.1]{Rein2007Review}).

    Therefore, let $0\le t \le s <T_{\max}$ and $z=(x,v) \in \supp\mathring{f}^- \cup \supp\mathring{f}^+$ be arbitrary. Using the fundamental theorem of calculus, we obtain
    \begin{align}
    \label{EST:V}
        \tfrac 12 |V^\pm(s,t,z)|^2 \le \tfrac 12 (R_v^0)^2 
            + C \int_t^s |V^\pm(\tau,t,z)|\, \|E_f(\tau)\|_{L^\infty(\mathbb R^{d_x})} \,\mathrm d\tau.
    \end{align}
    Note that the term involving the magnetic field vanishes since 
    $$V^\pm(s,t,z) \cdot \big(V^\pm(s,t,z) \otimes \overline{B}(s,X^\pm(s,t,z))\big) = 0 .$$
    Using \cref{lem:E}\ref{est:E:Linf} for the case $d_x=1$, and recalling that $f$ is constant along the measure-preserving characteristic flow, we deduce
    \begin{align*}
        \|E_f(\tau)\|_{L^\infty(\mathbb R^{d_x})} \le \|\rho_f(\tau)\|_{L^1(\mathbb R^{d_x})} 
        \le \|f(\tau)\|_{L^1(\mathbb R^{D})} = \|\mathring f\|_{L^1(\mathbb R^{D})}
    \end{align*}
    for all $\tau\in[0,T_{\max}[$.
    We thus obtain
    \begin{align}
    \label{EST:V2}
        \tfrac 12 |V^\pm(s,t,z)|^2 \le \tfrac 12 (R_v^0)^2 
            + C \int_t^s |V^\pm(\tau,t,z)| \,\mathrm d\tau.
    \end{align}
    Applying the quadratic version of Gronwall's lemma \cite[Theorem~5]{Dragomir}, we infer
    \begin{align}
    \label{EST:Z3}
        |V^\pm(s,t,z)| \le R_v^0 + C(s-t)           
    \end{align}
    Using the fundamental theorem of calculus, we further derive the estimate
    \begin{align}
    \label{EST:X}
        |X^\pm(s,t,z)| 
        \le R_x^0 + C \int_t^s |V(\tau,t,z)| \,\mathrm d\tau
        \le R_x^0 + C(s-t) + C(s-t)^2.
    \end{align}
    Combining \eqref{EST:Z3} and \eqref{EST:X}, we thus have
    \begin{align*}
        \XX(t) \le R_x^0 + Ct + Ct^2, \quad \VV(t) \le R_v^0 + Ct
    \end{align*}
    for all $t\in [0,T_{\max}[$. Thus, we conclude that $T_{\max}=\infty$.
    Moreover, by choosing $T$ sufficiently small, we ensure $\XX(t) < R_x$ and $\VV(t) < R_v$
    for all $t\in [0,T]$. As the constant $C$ depends only on $\mathring{f}$ and $d$, so does this choice of $T$.
    In particular, the solution $f$ depends only on $\overline{B}\vert_{[0,T]\times\overline{\Omega_x}} = B$ but not on the extension $\overline{B}$.
\end{proof}

\begin{rem}
    We point out that it is unclear whether the case $d=(3,3)$ can be handled in the same way. Solely following the proof of \cite[Theorem~13]{Knopf_Diss}, merely the estimate
    \begin{align}
        \label{EST:V33}
        \VV(t) \le C_B (1+t)^{21/2}
    \end{align}
    can be established with a constant $C_B$ depending on $\mathring f$, $R_x$, $R_v$ and $\|B\|_{L^2(0,T;L^\infty(\Omega_x))}$. Based on this inequality, the final time $T$ cannot be chosen independent of $\|B\|_{L^2(0,T;L^\infty(\Omega_x))}$. However, as we want to use magnetic fields that are defined on the interval $[0,T]$, where $T>0$ is supposed to be independent of the choice of $B$, this leads to a chicken and egg problem. It is an interesting open question whether an estimate in the fashion of \eqref{EST:V33} with a constant independent of $B$ can be established also in the case $d=(3,3)$.
\end{rem}

\begin{lemma}[Uniform bounds for classical solutions] \label{lem:f:bound}
    Let $K>0$ be arbitrary, let $B \in C^1([0,T]\times \overline{\Omega_x};\RR^{d_B})$ with $\|B\|_{W^{1,\infty}([0,T]\times\Omega_x)}\le K$, and let $f$ be the corresponding classical solution of \eqref{VP:PM:d} to the initial data $\mathring f$.
    Then, there exists a constant $c>0$ depending only on $\mathring{f}^\pm$, $R_x$, $R_v$, $d$ and $K$ such that 
    \begin{align*}
        \|f^\pm\|_{W^{1,\infty}([0,T]\times\RR^D)} \le c.
    \end{align*}
\end{lemma}

\begin{proof}
    Since $f$ is constant along the characteristic flow, we have 
    $\|f^\pm(t)\|_{L^\infty(\RR^D)} = \|\mathring{f}^\pm\|_{L^\infty(\RR^D)}$ for all $t\in [0,T]$.
    Moreover, using the estimates for the electric field presented in \cref{lem:E}\ref{est:E:Linf} and \ref{est:E:W1inf}, a uniform bound on $\|D_{(x,v)}f^\pm\|_{L^\infty([0,T]\times\RR^D)}$ 
    can be established by proceeding exactly as in \cite[Lemma~6]{Knopf2018OCP_VP}.
\end{proof}

\smallskip

\begin{lemma}[Local Lipschitz estimate for classical solutions] \label{lem:flip}
    Let $K>0$ be arbitrary. For $i\in\{1,2\}$, let $B_i \in C^1([0,T]\times \overline{\Omega_x};\RR^{d_x})$ with $\|B_i\|_{W^{1,\infty}([0,T]\times\Omega_x)}\le K$, and let $f_i$ denote the corresponding classical solution of \eqref{VP:PM:d} to the initial data $\mathring f$.
    Then, there exists a constant $L>0$ depending only on $R_x$, $R_V$, $\mathring{f}^\pm$ and $K$ such that
    \begin{align*}
        \|f_1^\pm - f_2^\pm\|_{C([0,T];L^2(\RR^D))} \le L\, \|B_1-B_2\|_{L^2([0,T]\times\RR^D)}.
    \end{align*}
\end{lemma}

\begin{proof}
    Let us write $\bar B\coloneqq  B_1-B_2$ and $\bar f^\pm = f_1^\pm - f_2^\pm$, and let $C>0$ denote a generic constant depending only on $R, \mathring{f}$ and $K$ that may change its value from line to line. Then $\bar f^\pm$ fulfills the Vlasov equations
    \begin{align*}
        \partial_t \bar f^\pm 
        + \mu_x^\pm v \cdot\partial_x \bar f^\pm 
        + \mu_v^\pm (E_{f_1} 
        + \mu_x^\pm v \otimes B_1) \cdot \partial_v \bar f^\pm 
        + \mu_v^\pm (E_{\bar f} + \mu_x^\pm v\otimes \bar B) \cdot \partial_v f_2^\pm
        =0,
    \end{align*}
    We now multiply these equations by $\bar f^\pm$, respectively. Integrating the resulting equations over $\RR^D$ and using integration by parts as well as \cref{lem:E}\ref{est:E:L2} and \cref{lem:f:bound}, we infer
    \begin{align*}
        &\frac{\mathrm d}{\mathrm dt} \|\bar f^\pm\|_{L^2(\RR^D)}^2
        = -2 \int_{\RR^D} 
            \big( \mu_v^\pm (E_{\bar f} + \mu_x^\pm v\otimes \bar B) \cdot \partial_v f_2^\pm \big)
            \bar f^\pm \mathrm dz 
        \\
        &\le C \big(\|E_{\bar f}\|_{L^2(\RR^D)} +
            \|\bar B\|_{L^2(\RR^D)}\big) \|\partial_v f_2^\pm\|_{L^\infty(\RR^D)} 
                \|\bar f^\pm\|_{L^2(\RR^D)}
        \\
        &\le C \|\bar f\|_{L^2(\RR^D)}^2
            + C \|\bar B\|_{L^2(\RR^D)} \|\bar f^\pm\|_{L^2(\RR^D)}
    \end{align*}
    on $[0,T]$. After integrating this estimate with respect to time from $0$ to $t$, we first use the standard version and then the quadratic version of Gronwall's lemma \cite[Theorem~5]{Dragomir} to deduce the claim.
\end{proof}

\subsection{Strong well-posedness for magnetic fields of class \texorpdfstring{$W^{1,\infty}$}{W\^{}\{1,\textbackslash infty\}}}
\label{sec:WellPosedness_strong}

We now prove strong well-posedness of the Vlasov--Poisson system 
\eqref{VP:PM:d} on $[0,T]$ 
for external magnetic fields that are Lipschitz continuous, i.e., of class $W^{1,\infty}$.

\begin{theorem}[Strong well-posedness for magnetic fields of class $W^{1,\infty}$]\label{thm:existence_strong_solutions}
	There exists $T>0$ depending only on $\mathring{f}^\pm$, $R_x$, $R_v$ and $d$ such that 
	for all $B \in W^{1,\infty}([0,T]\times\Omega_x;\RR^{d_B})$, 
	  there exists a unique strong solution $f=(f^-,f^+)$ with 
    $f^\pm \in W^{1,\infty}\left([0,T]\times\RR^D\right)$ 
    of the Vlasov--Poisson system \eqref{VP:PM:d} to the
    initial data $\mathring f$ with
	$\supp f^\pm(t) \subset \overline\Omega$ for all $t \in [0,T]$.
\end{theorem}

\begin{proof}
    Let $B\in W^{1,\infty}([0,T]\times\Omega_x;\RR^{d_B})$ be arbitrary. We now choose a sequence $(B_k)_{k\in\mathbb N} \subset C^1([0,T]\times \overline{\Omega_x};\RR^{d_B})$ with 
    \begin{align*}
        &\|B_k\|_{W^{1,\infty}([0,T]\times \Omega_x)} \le \|B\|_{W^{1,\infty}([0,T]\times \Omega_x)}
        \quad\text{for all $k\in\mathbb N$},\\
        &\text{and}\quad 
        B_k\to B \quad\text{in $L^2([0,T]\times \overline{\Omega_x})$ as $k\to\infty$.}
    \end{align*}
    According to \cref{prop:existence_classical_solutions}, for every $k\in\mathbb N$, there exists a unique classical solution $f_k=(f_k^-,f_k^+)$ of the Vlasov--Poisson system \eqref{VP:PM:d} on the interval $[0,T]$ (with $T$ independent of $k$) to the
    initial data $\mathring f$ with
	$\supp f_k^\pm(t) \subset \overline\Omega$ for all $t \in [0,T]$.
    We further know from \cref{lem:f:bound} that the sequences $(f_k^\pm)_{k\in\mathbb N}$ are bounded in the space $W^{1,\infty}\left([0,T]\times\RR^D\right)$. Hence, there exists $f=(f^-,f^+)$ with $f \in W^{1,\infty}\left([0,T]\times\RR^D\right)$ such that
    \begin{align*}
        f_k^\pm \overset{\ast}{\rightharpoonup} f^\pm, \quad
        \partial_t f_k^\pm \overset{\ast}{\rightharpoonup} \partial_t f^\pm, \quad
        \partial_{z_i} f_k^\pm \overset{\ast}{\rightharpoonup} \partial_{z_i} f^\pm
        \quad\text{in $L^\infty\big([0,T]\times \RR^D\big)$},
    \end{align*}
    for all $i\in\{1,...,D\}$, as $k\to\infty$ along a non-relabeled subsequence. In particular, since $\supp f_k^\pm(t) \subset \overline\Omega$ for all $t \in [0,T]$, this implies that $\supp f^\pm(t) \subset \overline\Omega$ for all $t \in [0,T]$. Let now $r>0$ such that $\overline\Omega\subset B_r(0)$.
    Furthermore, it follows from \cref{lem:f:bound} that the sequences $(f_k^\pm)_{k\in\mathbb N}$ are also bounded in $W^{1,8}\big([0,T]\times B_r(0)\big)$.
     Since $8>D+1$, the embedding $W^{1,8}([0,T]\times B_r(0)) \hookrightarrow C([0,T]\times B_r(0))$ is compact. 
     As the supports of $f^\pm$ and all $f_k^\pm$ remain contained in $\overline\Omega\subset B_r(0)$, we thus have
    \begin{align*}
        f_k^\pm \to f^\pm \quad\text{in $C([0,T]\times \RR^D)$}
    \end{align*}
    as $k\to\infty$ up to subsequence extraction. In particular, this proves that $f$ satisfies the initial condition.
    Using \cref{lem:E}\ref{est:E:Linf} and \ref{est:E:W1inf}, we further obtain
    \begin{align*}
        E_{f_k} \to E_f \quad\text{in  $C([0,T]\times \overline \Omega_x;\RR^{d_x})$}.
    \end{align*}
    Multiplying the Vlasov equation from \eqref{VP:PM:d} written for $B_k$ and $f_k$ with an arbitrary test function $\varphi\in C^\infty_c([0,T]\times \RR^D)$ and integrating the resulting equation over $[0,T]\times \RR^D$, we can use the above convergences to pass to the limit $k\to\infty$. By means of the fundamental lemma of the calculus of variations, we conclude that $f \in C([0,T]\times\RR^D) \cap W^{1,\infty}([0,T]\times\RR^D)$ satisfies the Vlasov equation \eqref{VP:PM:d} almost everywhere in $[0,T]\times \RR^D$. This means that $f$ is a strong solution of the Vlasov--Poisson system to the magnetic field $B$ and the initial data $\mathring{f}$, which satisfies the desired support condition. 

    The uniqueness of this strong solution follows by a straightforward Gronwall argument similar to the one in the proof of \cref{lem:flip}. Thus, the proof is complete.
\end{proof}

\subsection{The control-to-state operator and its properties}
\label{sec:ControlToStateOperator}

\begin{definition}[Control-to-state operator]
\label{def:ControlToState}
    For any magnetic field $B\in W^{1,\infty}([0,T]\times\Omega_x;\RR^{d_B})$, let $f_B=(f_B^-,f_B^+)$ with $f_B^\pm \in W^{1,\infty}([0,T]\times\RR^D)$ be the unique strong solution of the Vlasov--Poisson system \eqref{VP:PM:d}. Then, the operator 
    \begin{align*}
        f. = (f.^-,f.^+) : W^{1,\infty}([0,T]\times\Omega_x;\RR^{d_B}) \to W^{1,\infty}([0,T]\times\Omega)^2,
        \quad B \mapsto f_B
    \end{align*}
    is called the \emph{control-to-state operator}.
\end{definition}

We shall now derive important properties of the control-to-state operator. In particular, due to a lack of higher regularity, we prove local Lipschitz continuity and Fréchet differentiability with respect to weaker norms than the $W^{1,\infty}$-norm. To this end, recall that $W^{1,\infty}([0,T]\times\Omega)\subset C([0,T];L^2(\RR^D))\subset C([0,T];H^{-1}(\RR^D))$ in the sense of continuous embeddings (and by extension by zero outside of $\Omega$).

\begin{theorem}[Boundedness and local Lipschitz continuity]\label{thm:cso:lipschitz}
    Let $K>0$ be arbitrary. Then, there exist constants $C, L \ge 0$ depending only on $\mathring{f}$, $R_x$, $R_v$, $d$ and $K$ such that: 
    \begin{enumerate}[label=\textnormal{(\alph*)}]
        \item For all $B\in W^{1,\infty}([0,T]\times\Omega_x;\RR^{d_B})$ with $\|B\|_{W^{1,\infty}([0,T]\times\Omega_x)}\le K$,
        \begin{align*}
            \|f_B^\pm\|_{W^{1,\infty}([0,T]\times\Omega)} \le C.
        \end{align*}
        \item For all $B_1,B_2\in W^{1,\infty}([0,T]\times\Omega_x;\RR^{d_B})$ with $\|B_i\|_{W^{1,\infty}([0,T]\times\Omega_x)}\le K$ for $i=1,2$,
        \begin{align*}
            \| f_{B_1}^\pm - f_{B_2}^\pm\|_{C([0,T];L^2(\RR^D))} \le L \|B_1 - B_2\|_{L^2([0,T]\times\Omega_x)}.
        \end{align*}
    \end{enumerate}
\end{theorem}

\begin{proof}
    As in the proof of \cref{thm:existence_strong_solutions}, we approximate any $B\in W^{1,\infty}([0,T]\times\Omega_x;\RR^{d_B})$ by a sequence $(B_k)_{k\in\mathbb N} \subset C^1([0,T]\times \overline{\Omega_x};\RR^{d_B})$ which satisfies 
    \begin{align*}
        &\|B_k\|_{W^{1,\infty}([0,T]\times \Omega_x)} \le \|B\|_{W^{1,\infty}([0,T]\times \Omega_x)} \le K
        \quad\text{for all $k\in\mathbb N$},\\
        &\text{and}\quad 
        B_k\to B \quad\text{in $L^2([0,T]\times \overline{\Omega_x})$ as $k\to\infty$.}
    \end{align*}
    Then, as in the proof of \cref{thm:existence_strong_solutions}, the corresponding classical solutions $f_{B_k}$ satisfy 
    \begin{align*}
        & \|f_{B_k}^\pm\|_{W^{1,\infty}([0,T]\times\Omega)} \le c, \\
        & f_{B_k}^\pm \overset{\ast}{\rightharpoonup} f_B^\pm, \quad
        \partial_t f_{B_k}^\pm \overset{\ast}{\rightharpoonup} \partial_t f_B^\pm, \quad
        \partial_{z_i} f_{B_k}^\pm \overset{\ast}{\rightharpoonup} \partial_{z_i} f_B^\pm
        \quad\text{in $L^\infty\big([0,T]\times \RR^D\big)$}, \\
        &\text{and}\quad 
        f_{B_k}^\pm \to f_B^\pm \quad\text{in $C([0,T]\times \RR^D)$}
    \end{align*}
    after subsequence extraction. Hence, (a) follows from \cref{lem:f:bound} and the weak-$^*$ lower semicontinuity of the $L^\infty([0,T]\times\Omega)$-norm, whereas (b) follows from \cref{lem:flip} by passing to the limit $k\to\infty$.
\end{proof}

\begin{theorem}[Frech\'et differentiability]\label{thm:cso:frechet}
    The control-to-state operator is Fr\'echet differentiable in the following sense: For every $B\in W^{1,\infty}([0,T]\times\Omega_x;\RR^{d_B})$ there exists a linear operator
    \begin{align*}
        \left\{
        \begin{aligned}
        f_B'[\cdot]
        = \big(f_B^{\prime-}[\cdot],f_B^{\prime+}[\cdot]\big):
        W^{1,\infty}([0,T]\times\Omega_x;\RR^{d_B}) &\to W^{1,\infty}([0,T]\times\RR^D)^2,
        \\
        H &\mapsto f_B'[H]
        \end{aligned}
        \right.
    \end{align*}
    such that
    \begin{align*}
        \frac{\|f_{B+H}^\pm - f_B^\pm - f^{\prime\pm}_B[H]\|_{C([0,T];H^{-1}(\RR^D))}}
        {\|H\|_{L^2([0,T]\times\Omega_x)}}
        \to 0
        \quad\text{as $\|H\|_{L^2([0,T]\times\Omega_x)} \to 0$.}
    \end{align*}
\end{theorem}

\medskip

\begin{proof}
    Let $B,H\in W^{1,\infty}([0,T]\times\Omega_x;\RR^{d_B})$ be arbitrary. Let $f_B'[H] \in H^1([0,T]\times \RR^D)$ be the unique strong solution of the system 
    \begin{align}
        \label{VP:lin}
        \left\{
        \begin{aligned}
            &\partial_t f^\pm + \mu_x^\pm v \cdot\partial_x f^\pm 
            + \mu_v^\pm (E_{f_B} + \mu_x^\pm v\otimes B) \cdot \partial_v f^\pm
            =
            - \mu_v^\pm (E_f + \mu_x^\pm v\otimes H) \cdot \partial_v f_B^\pm,
            \\
            & f^\pm\vert_{t=0} = 0,
        \end{aligned}
        \right.
    \end{align}
    and let $f_{\mathcal R}[H]\in H^1([0,T]\times \RR^D)$ be the unique strong solution of the system  
    \begin{align}
        \label{VP:rem}
        \left\{
        \begin{aligned}
            &\partial_t f^\pm + \mu_x^\pm v \cdot\partial_x f^\pm 
            + \mu_v^\pm (E_{f_B} + \mu_x^\pm v\otimes B) \cdot \partial_v f^\pm
            =
            - \mu_v^\pm E_f \cdot \partial_v f_B^\pm 
            + \mathcal R,
            \\
            & f^\pm\vert_{t=0} = 0,
        \end{aligned}
        \right.
    \end{align}
    where 
    \begin{equation*}
        \mathcal R \coloneqq  -\mu_v^\pm (E_{f_{B+H} - f_B} 
          - \mu_x^\pm v\otimes H) \cdot 
            \big( \partial_v f_{B+H}^\pm - \partial_v f_{B}^\pm \big).
    \end{equation*}
    We point out that the existence and the uniqueness of these strong solutions $f_B'[H]$ and $f_{\mathcal R}[H]$ can be established exactly as in \cite[Corollary~2]{Knopf2018OCP_VP}, where the strong well-posedness of such linear Vlasov equations was investigated in the case $d=(3,3)$. To obtain more information about the supports of $f_B^{\prime\pm}[H](t)$ for all $t\in [0,T]$, we define $Z_B^\pm=(X_B^\pm,V_B^\pm)(s,t,z)$ to denote the unique solutions of the characteristic systems
    \begin{align}
        \label{charsys:B}
        \dot x(s) = v(s),\quad \dot v(s) =  \mu_v^\pm E_{f_B}\big(s,x(s)\big) + \mu_x^\pm \mu_v^\pm v(s) \otimes B\big(s,x(s)\big)
    \end{align}
    satisfying the initial conditions $Z_B^\pm(t,t,z) = z$, respectively, for any given $z=(x,v)\in\mathbb R^6$. 
    Following the proof of \cite[Corollary~2]{Knopf2018OCP_VP}, we deduce that $f_B^{\prime\pm}[H]$ can be represented as
    \begin{align*}
        f_B^{\prime\pm}[H](t,z)
        = - \int_0^t \big[\mu_v^\pm (E_f + \mu_x^\pm v\otimes H) \cdot \partial_v f_B^\pm\big]\big(s,Z_B^\pm(s,t,z)\big) \,\mathrm ds
    \end{align*}
    for all $(t,z)\in [0,T]\times\RR^D$. Recalling that $f_B$ is constant along the characteristic flow $Z_B$, we conclude 
    \begin{align*}
        \supp f_B^{\prime\pm}[H](t)
        \subseteq \bigcup_{s\in [0,T]} \supp f_B^{\pm}\big(s,Z_B(s,t,\cdot)\big)
        = \supp f_B^{\pm}(t)
        \subseteq \overline\Omega
    \end{align*}
    for all $t\in[0,T]$. Proceeding analogously, we also obtain 
    $\supp f^\pm_{\mathcal R}[H](t) \subset \overline\Omega$ for all $t\in[0,T]$.   
    Furthermore, is easy to check that $f_{B+H}^\pm - f_B^\pm - f^{\prime\pm}_B[H]$ is a strong solution of \eqref{VP:rem} and thus, because of uniqueness, 
    \begin{align}
        \label{ID:FR}
        f_{B+H}^\pm - f_B^\pm - f^{\prime\pm}_B[H] = f_{\mathcal R}^\pm[H] .
    \end{align}
    It remains to derive a suitable bound on $f_{\mathcal R}[H]$. 
    To estimate $\|f_{\mathcal R}[H]\|_{C([0,T];H^{-1}(\RR^D))}$, we use a similar idea as in \cite{Robert}, where such techniques were used in a uniqueness proof for weak solutions. 
    For any $t\in [0,T]$, let $\phi(t,\cdot) = (\phi^-(t,\cdot),\phi^+(t,\cdot)) \in H^1(\RR^D)^2$ denote the unique solution of Poisson's equation
    \begin{align*}
        \left\{
        \begin{aligned}
            - \Delta_z \phi^\pm(t,\cdot) &= f^\pm_{\mathcal R}[H](t,\cdot)
            &&\quad\text{in $\RR^D$},
            \\
            \underset{|z|\to\infty}{\lim}\; \phi^\pm(t,z) &= 0 ,
        \end{aligned}
        \right.
    \end{align*}
    where $\Delta_z$ denotes the Laplace operator in $\RR^D$. From the regularity of $f_{\mathcal R}[H]$, we infer that $\phi \in H^1(0,T;H^2(\Omega)) \cap L^2(0,T;H^3(\Omega))$ by means of elliptic regularity theory. Let now $\varphi \in H^1(\RR^D)$ and $t\in [0,T]$ be arbitrary. Using integration by parts, we deduce
    \begin{align*}
        &\big| \big(f^\pm_{\mathcal R}[H](t), \varphi\big)_{L^2(\RR^D)}\big|
        = \big|\big(-\Delta_z \phi^\pm(t), \varphi\big)_{L^2(\RR^D)}\big|
        \\
        &\quad= \big|\big(\nabla_z \phi^\pm(t), \nabla_z \varphi\big)_{L^2(\RR^D)}\big|
        \le \|\nabla_z \phi^\pm(t)\|_{L^2(\RR^D)} \|\varphi\|_{H^1(\RR^D)}.
    \end{align*}
    Here, $\nabla_z$ denotes the gradient operator in $\RR^D$. Note that the integration by parts does not produce any boundary terms as the solution $\phi^\pm$ decays sufficiently fast as $|z|\to\infty$.
    Taking the supremum over all $\varphi \in H^1(\RR^D)$ with $\|\varphi\|_{H^1(\Omega)} \le 1$ and all $t\in [0,T]$, we conclude
    \begin{align}
    \label{EST:H-1}
        \|f^\pm_{\mathcal R}[H]\|_{C([0,T];H^{-1}(\RR^D))}
        \le \|\nabla_z \phi^\pm\|_{C([0,T];L^2(\RR^D))}.
    \end{align}
    It thus remains to derive a suitable bound on the right-hand side of this inequality. 
    Notice that $f^\pm_{\mathcal R}[H](0,\cdot) = 0$ implies $\phi^\pm(0,\cdot) = 0$.
    Using $\phi^\pm$, the system \eqref{VP:rem} can be rewritten as
    \begin{align*}
        \label{VP:phi}
        \left\{
        \begin{aligned}
            &- \partial_t \Delta_z \phi^\pm 
            - \mu_x^\pm v \cdot\partial_x \Delta_z \phi^\pm  
            - \mu_v^\pm (E_{f_B} + \mu_x^\pm v\otimes B) \cdot \partial_v \Delta_z \phi^\pm  =
            \mu_v^\pm E_{\Delta_z \phi} \cdot \partial_v f_B^\pm 
            + \mathcal R,
            \\
            & \phi^\pm \vert_{t=0} = 0,
        \end{aligned}
        \right.
    \end{align*}
    Multiplying by $\phi^\pm$ and using the fundamental theorem of calculus, we obtain
    \begin{align*}
        \frac 12 \|\nabla_z \phi^\pm(t)\|_{L^2(\RR^D)}^2
        = \int_0^t \int_{\RR^D} & \big\{ \mu_x^\pm v \cdot\partial_x \Delta_z \phi^\pm \, \phi^\pm
            + \mu_v^\pm (E_{f_B} + \mu_x^\pm v\otimes B) \cdot \partial_v \Delta_z \phi^\pm \, \phi^\pm \\[-1ex]
            &
            + \mu_v^\pm E_{\Delta_z\phi} \cdot \partial_v f_B^\pm \, \phi^\pm+ \mathcal R \, \phi^\pm \big\} \,\mathrm dz \mathrm ds
    \end{align*}
    Via integration by parts, this equation can be reformulated as 
    \begin{align}
        \label{EST:DPHI}
        \frac 12 \|\nabla_z \phi^\pm(t)\|_{L^2(\RR^D)}^2 
        = &\int_0^t \int_{\RR^D} \Big\{ \mu_x^\pm \partial_x \phi^\pm \cdot \partial_v \phi^\pm
            - \mu_v^\pm f_B^\pm E_{\Delta_z\phi} \cdot \partial_v \phi^\pm  
            \\[1ex]\notag
            & + \mu_v^\pm \big[ D_x (E_{f_B} + \mu_x^\pm v\otimes B) \partial_x \phi^\pm \big] \cdot \partial_v \phi^\pm
            \\[1ex]\notag
            &+ \mu_v^\pm\mathcal (f_{B+H} - f_B ) \big[ E_{f_{B+H}-f_B} + \mu_x^\pm v\otimes H \big] \, \cdot \partial_v \phi^\pm \Big\} \,\mathrm dz \mathrm ds
    \end{align}
    for all $t\in[0,T]$. Due to Gau\ss's theorem and Fubini's theorem, we have
    \begin{align*}
        \rho_{\Delta_v\phi}(t,x) 
        = 0
        \quad\text{and}\quad
        \int_{\RR^{d_x}} \rho_{\Delta_x\phi}(t,x) \,\mathrm dx 
        = 0 
    \end{align*}
    for all $t\in[0,T]$ and $x\in\RR^{d_x}$.
    Since $\Delta_z \phi = \Delta_x\phi + \Delta_v\phi$, this implies
    \begin{align*}
        E_{\Delta_z\phi} = E_{\Delta_x \phi} 
        = - \partial_x (-\Delta_x)^{-1} \int_{\RR^{d_v}} \Delta_x (\phi^+ - \phi^-) \,\mathrm dv  
        = \int_{\RR^{d_v}} \partial_x (\phi^+ - \phi^-) \,\mathrm dv 
    \end{align*}
    Following the line of argument in \cite{Robert} and meanwhile using the established uniform control on the support of $f_{\mathcal R}[H]$, we infer that
    \begin{align}
        \label{EST:E:CZ}
        \|E_{\Delta_z \phi}\|_{L^2(\Omega)} \le C
        \sum_\pm \|\partial_x \phi^\pm\|_{L^2(\RR^D)}.
    \end{align}
    Using Hölder's inequality along with \eqref{EST:E:CZ}, \cref{lem:f:bound} and \cref{lem:flip} to bound the right-hand side of \eqref{EST:DPHI}, we deduce the estimate
    \begin{align}
        \label{EST:DPHI:2}
        &\sum_\pm \|\nabla_z \phi^\pm(t)\|_{L^2(\RR^D)}^2 
        \\\notag
        &\le C \sum_\pm \int_0^t \big(\|\nabla_z \phi^\pm(s)\|_{L^2(\RR^D)}^2 
        + \|H\|_{L^2([0,T]\times\RR^{d_x})}^2 \|\nabla_z \phi^\pm(s)\|_{L^2(\RR^D)}\big)\,\mathrm ds,
    \end{align}
    where $C>0$ does not depend on $H$.
    Applying first the standard version and then the quadratic version of Gronwall's lemma \cite[Theorem~5]{Dragomir}, we conclude
    \begin{align*}
        \sum_\pm \|\nabla_z \phi^\pm(t)\|_{L^2(\RR^D)} \le C' \|H\|_{L^2([0,T]\times\RR^{d_x})}^2,
    \end{align*}
    where $C'>0$ does not depend on $H$.
    In view of \eqref{ID:FR} and \eqref{EST:H-1}, we obtain
    \begin{align*}
        \|f_{B+H}^\pm - f_B^\pm - f^{\prime\pm}_B[H]\|_{C([0,T];H^{-1}(\RR^D))} \le C' \|H\|_{L^2([0,T]\times\RR^{d_x})}^2,
    \end{align*}
    which proves the claim.
\end{proof}

\subsection{The set of admissible controls}
\label{sec:AdmissibleControls}
We now introduce the set of magnetic fields that are admissible for our optimal control problem. Let $\mathcal{V}$ denote the Hilbert space $H^1([0,T]\times \Omega_x; \RR^{d_B})$ on which we define the weighted scalar product
\begin{align}
    \label{def:Vscp}
    (B_1,B_2)_{\mathcal V} &\coloneqq  (B_1,B_2)_{L^2([0,T]\times\Omega_x)} 
        + \kappa_t (\partial_t B_1 , \partial_t B_2)_{L^2([0,T]\times\Omega_x)} 
    \\ \notag
    &\qquad + \kappa_x (D_x B_1,D_x B_2)_{L^2([0,T]\times\Omega_x)}
\end{align}
for a given vector $\kappa = (\kappa_t,\kappa_x)$ with $\kappa_t,\,\kappa_x>0$ and all $B_1,B_2\in\mathcal V$. The induced norm is given by 
\begin{align}
    \label{def:Vnorm}
    \|B\|_{\mathcal V} \coloneqq  (B,B)_{\mathcal V}^{\frac 12} \quad\text{for all $B\in\mathcal V$}.
\end{align}

\begin{definition}
    For every $i\in\{0,...,d_x\}$, $j\in \{1,...,d_B\}$, let $a_{ij},b_{ij} \in L^\infty([0,T]\times \Omega_x)$ with $a_{ij} \le 0 \le b_{ij}$ a.e.~in $[0,T]\times \Omega_x$ be given functions.
	Then the set
	\begin{align*}
		\Bad \coloneqq \left\{ B \in \mathcal{V}
		\;\middle|\;
        \begin{aligned}
        &a_{0j} \le \partial_t B_j \le b_{0j}   
        \;\text{a.e.~in $[0,T]\times \Omega_x$},\\
        &a_{ij} \le \partial_{x_i} B_j \le b_{ij}   
        \;\text{a.e.~in $[0,T]\times \Omega_x$}\\
        &\text{for all $i\in\{1,...,d_x\}$ and $j\in \{1,...,d_B\}$}
        \end{aligned}
        \right\}
	\end{align*}
	is called the \emph{set of admissible controls}.
	\label{def:admissible_controls}
\end{definition}
On the one hand, from a physical point of view, the choice of the set of admissible controls is motivated by the fact that realistic external magnetic fields will not have arbitrarily large spatial gradients and do not change arbitrarily fast in time. In particular, as illustrated in \cite[Fig.~5]{Bartsch2021Numerical}, choosing $L^2([0,T]\times \Omega_x; \RR^{d_B})$ instead of $\mathcal V$ might lead to discontinuous (locally) optimal controls, which are not physically reasonable. On the other hand, from a mathematical point of view, this choice is useful as the set $\Bad$ has the following properties:

\begin{lemma}\label{lem:bad}   
    There exists a constant $C>0$ depending only on $R_x$, $d$, $\kappa$ and the box constraints $a_{ij}$ and $b_{ij}$, $i\in\{0,...,d_x\}$, $j\in \{1,...,d_B\}$, such that
    \begin{align*}
        \|B\|_{W^{1,\infty}([0,T]\times\Omega_x)} 
        \le C\big(1 + \|B\|_{\mathcal V}\big)
        \quad\text{for all $B\in \Bad$.}
    \end{align*}
    In particular, this means that every field $B\in \Bad$ is Lipschitz continuous on $[0,T]\times\Omega_x$.
\end{lemma}

\begin{proof}
For all $B\in \Bad \cap C^1([0,T]\times \overline{\Omega_x};\RR^{d_B})$, $(t,x)\in [0,T]\times\Omega_x$ and $j\in \{1,...,d_B\}$,
\begin{align*}
    |B_j(t,x)| 
    &\le |\Omega_x|^{-1}  \int_{\Omega_x} \left(|B_j(t,x) - B_j(t,y)| + |B_j(t,y)|\right) \,\mathrm dy \\
    & \le 2R_x\, \|\nabla_x B_j(t)\|_{L^\infty(\Omega_x)} + |\Omega_x|^{-\frac 12}\|B_j\|_{C([0,T];L^2(\Omega_x))} \le C\big(1 + \|B\|_{\mathcal V}\big)
\end{align*}
for some constant $C>0$ depending only on $R_x$, $d$, $\kappa$ and the box constraints $a_{ij}$ and $b_{ij}$. Since $C^1([0,T]\times \overline{\Omega_x};\RR^{d_B}) \subset \mathcal V$ is dense, we conclude
\begin{align*}
    \|B\|_{L^\infty([0,T]\times \Omega_x)} \le C+C\|B\|_{\mathcal V}
\end{align*}
for all $B\in\Bad$. 
Along with the box constraints incorporated in $B\in \Bad$, this proves the claim.
\end{proof}

\subsection{Formulation of the optimal control problem}
\label{sec:Formulation_OPC}

Now, we discuss the formulation of our cost functional that includes penalization terms for the distribution functions and for the control. 
As in \cite{Bartsch2019Theoretical,Bartsch2021MOCOKI}, we choose our cost functional in the style of ensemble control problems; see, e.g., \cite{Brockett2012NotesLiouvilleControl} and the references therein. 
We define
\begin{align}
    J(f,B) &\coloneqq \sum_\pm 
    \int_0^T \int_{\Omega_x} \int_{\Omega_v}   \theta^\pm(t,x,v) f^\pm(t,x,v)\,\mathrm dv\mathrm dx \mathrm dt
    \label{eq:objective_function}\\ &\phantom{\coloneqq\;}+
    \sum_\pm \int_{\Omega_x} \int_{\Omega_v}  \varphi^\pm(x,v) f^\pm(T,x,v)\,\mathrm dv\mathrm dx +
    \frac{\alpha}{2} \|B\|^2_{\mathcal{V}}
    ,\notag
\end{align}
where $\alpha>0$ acts as a weight for the term $\|B\|^2_{\mathcal{V}}$, which penalizes the influence of the control.
The summands involving $f^\pm$ can be interpreted as weighted expected values which are often used in the context of kinetic equations in statistical mechanics.

We refer to the first summand in \eqref{eq:objective_function} as the tracking term (also called integrated cost)
and to the second summand as the final observation 
(also called terminal cost). The purpose of the tracking term is to penalize deviations of the particle distribution from a prescribed (time-dependent) mean phase-space profile $z_{\mathfrak{d}}^\pm(t)$. Similarly, the final observation term penalizes deviations from a prescribed mean phase-space profile $z_{T}^\pm$ at the final time $T$.
A reasonable choice would be $\theta^\pm(t,z) = \Theta^\pm \left( | z-z_\mathfrak{d}^\pm(t)|^2 \right)$, where $\Theta^\pm$ is continuously differentiable and monotonically increasing.
Then the global minimum of the tracking part would be achieved when all particles have position and velocity equal to $z_\mathfrak{d}^\pm(t)$ at every time $t\in[0,T]$. 
Similarly, the final observation term could be chosen as $\varphi^\pm(z) = \Phi^\pm\left(|z - z_T^\pm|^2\right)$ with a continuously differentiable, increasing function $\Phi^\pm$ to achieve that the mean phase-space configuration of the particles at final time is close to $z_T^\pm$.
However, to keep $\theta^\pm$ and $\varphi^\pm$ as general as possible, we make the following assumption:

\begin{assumption}
	\label{asmpt:potentials}  
	We suppose that $\theta^\pm\in C^1\big([0,T]\times\Omega\big)$ and $\varphi^\pm\in C^1(\Omega)$
    are bounded from below, attaining for each $t \in [0,T]$ their global minimum at $z_\mathfrak{d}^\pm(t)\in\Omega$ and $z_T^\pm\in \Omega$.
    Moreover, they are assumed to be convex in a neighborhood of $z_\mathfrak{d}^\pm(t)$ and $z_T^\pm$, respectively. 
\end{assumption}

We point out that the local convexity assumption in a neighborhood around the minimizers is not necessary for our mathematical analysis.
However, it makes sense regarding the goal we want to achieve by our optimization, and it will also lead to qualitatively better numerical results.
Concerning the third term in \eqref{eq:objective_function}, 
we choose $\mathcal{V}$ as specified in Section \ref{sec:AdmissibleControls}. In particular, the norm $\|\cdot\|_{\mathcal V}$ is given by \eqref{def:Vnorm}.

Now, we can formulate our kinetic optimal control problem:
\begin{align}
		\min	J(f,B)
		\qquad
		\text{ s.t. } 
		\qquad 
		(f, B) \in W^{1,\infty}\bigl([0,T]\times\mathbb{R}^D \bigr)^2 \times \Bad  
		\text{ satisfies } \eqref{VP:PM:d},
	\label{eq:optimization_problem}
\end{align}
where $J$ is defined in \eqref{eq:objective_function} and $\Bad$ is given by Definition \ref{def:admissible_controls}.

\subsection{Existence of an optimal control}

By means of the control-to-state operator from \cref{def:ControlToState}, we reformulate the optimal control problem \eqref{eq:optimization_problem} as the so-called \textit{reduced optimal control problem}:
\begin{align}
	\min_{B \in \Bad} \widetilde{J}(B).
	\label{eq:reduced_optim_problem}
\end{align}
Here, $\widetilde{J}(B) \coloneqq J(f_B,B)$ defines the so-called \textit{reduced cost functional}.
We now prove the existence of solutions to \eqref{eq:reduced_optim_problem}:

\begin{theorem}[Existence of optimal solutions]
	Let $\theta^\pm$ and $\varphi^\pm$ fulfill \cref{asmpt:potentials}.
	The optimal control problem \eqref{eq:reduced_optim_problem} possesses at least one globally optimal solution $\bar{B}$.
\end{theorem}

\begin{proof}
    Using that $\|f_B^\pm(t)\|_{L^1(\RR^D)} = \|\mathring{f}^\pm\|_{L^1(\RR^D)}$ for all $t\in [0,T]$, we derive the estimate 
    \[ 
        \widetilde J(B) \ge \sum_\pm(T\inf\theta^\pm+\inf\varphi^\pm)\|\mathring f^\pm\|_{L^1(\RR^D)} + \frac{\alpha}{2} \|B\|^2_{\mathcal{V}}
        \quad\text{for all $B\in\Bad$}.
    \]
    Hence, the functional $\widetilde J$ is coercive and bounded from below. Thus, any minimizing sequence $\left(B_k\right)_{k \in \mathbb{N}}\subset \Bad$ of $\widetilde J$ is bounded in $\mathcal{V}$, and therefore also in $W^{1,\infty}([0,T]\times\Omega_x)$ by virtue of \cref{lem:bad}. Thus, such a sequence converges, after extracting a suitable subsequence, weakly in $\mathcal{V}$ to some $\bar B$, while the partial derivatives of $B_k$ converge weakly-* in $L^\infty([0,T]\times\Omega_x)$ to the partial derivatives of $\bar B$. In particular, $\bar B\in\Bad$ and $B_k\to \bar B$ strongly in $L^2([0,T]\times\Omega_x)$. Arguing as in the proof of \cref{thm:cso:lipschitz}, it follows that $f_{B_k}\to f_{\bar B}$ strongly in $C([0,T];L^2(\RR^D))$, and that $f_{\bar B}\in W^{1,\infty}([0,T];L^2(\RR^D))$. This and the uniform control on the support of $f_{B_k}^\pm(t)$ (cf.~\cref{thm:existence_strong_solutions}), directly imply that the terms in $\widetilde{J}$ containing $f_{B_k}$ converge to the respective expressions with $f_{\bar B}$ instead. Hence, as the $\mathcal V$-norm is weakly lower semicontinuous, we conclude
    \[
        \widetilde J(\bar B)\le\liminf_{k\to\infty} \widetilde J(B_k).
    \]
    This proves that $\bar B$ is a minimizer of $\widetilde{J}$.
\end{proof}

We remark that in our kinetic model \eqref{VP:PM:d} the Lorentz force $\mu_v^\pm (E_f + \mu_x^\pm v\otimes B)$ is multiplied by $\partial_v f^\pm$, and this bilinear structure makes our optimization problem non-convex and nonlinear. For this reason, in general, it is not possible to establish uniqueness of the minimizer $\bar B$. However, if certain parameters (here $T/\alpha$) are sufficiently small, then one can typically derive uniqueness. Such a result can also be established in our case, but we omit the proof and refer the interested reader to \cite{KnopfWeber2020_OCPVP_Coils} for details in a very similar setting.

It is further important to point out that in the above proof, we crucially needed the partial derivatives of the magnetic fields along a minimizing sequence to be a priori uniformly bounded in order to get a Lipschitz magnetic field in the limit, and in order to be able to pass to the limit in the Vlasov--Poisson system. In principle, one can imagine two ways to ensure this a priori property: 
First, our choice, namely to work with a quite weak norm for the magnetic fields that does not control the Lipschitz norm, and instead to impose box constraints on the partial derivatives of the magnetic field. Second, to impose no box constraints at all, but to choose a stronger norm for the magnetic fields that actually does control the Lipschitz norm. 
The advantage of the second approach, in view of what we will do below, would be that one would get an elliptic \textit{equation} as the first-order necessary optimality condition, from which one could bootstrap the regularity of an optimal magnetic field by classical elliptic regularity theory. A disadvantage would be that a more restrictive control space might lead to optimal controls that are quantitatively not as good as the ones from a more general control space. 
For instance, if higher order derivatives are punished by the cost functional, it may often be the most favorable option for the control to be close to constant. However, this might not be the best option to actually achieve the optimization goal.
Furthermore, the elliptic equation resulting from the first-order necessary optimality condition would be of order higher than two (since a stronger norm than $H^1$ is chosen for the magnetic fields in the first place), so that the second approach is more inconvenient from a numerical point of view. On the other hand, the first approach has the disadvantage of necessarily requiring box constraints and therefore, in general, only leading to a variational \textit{inequality}. 
However, it is, as we will discuss later, reasonable to ignore these box constraints in the numerical part, thus leading to an elliptic problem of order only two as the first-order necessary optimality condition. Therefore, from a numerical point of view, this first approach is the more attractive one, which is why we chose to follow it.

\subsection{First-order necessary optimality conditions}
\label{sec:OptimalitySystem}

Since the control-to-state operator $B\mapsto f_B$ is Fr\'echet differentiable (see \cref{thm:cso:frechet}) and the reduced cost functional $\widetilde J$ depends linearly and continuously on $f_B$, it follows by means of the chain rule that $\widetilde J$ is also Fr\'echet differentiable.

Let now $B\in\Bad$ be arbitrary.
To derive a formula for $\widetilde{J}'(B)$, which allows its numerical computation, as well as a first-order necessary optimality condition, we employ the so-called \textit{adjoint approach}.
Therefore, let us fix a cutoff function $\chi \in C^\infty_c(\RR^D)$ with $\chi = 1$ in $B_{R^*}(0)$, where $R^*>0$ is chosen sufficiently large, at least such that $\overline\Omega\subset B_{R^*}(0)$. Then the \textit{adjoint system}, which can be derived by the \textit{formal Lagrangian method}, reads as:%
\begin{subequations}\label{eq:adjoint_model}
	\begin{align}
		&\partial_t \lambda^\pm + \mu_x^\pm v\cdot\partial_x \lambda^\pm +\mu_v^\pm (E_{f_B}+\mu_x^\pm v\otimes B)\cdot\partial_v \lambda^\pm 
        = \Phi_{f_B,\lambda}^\pm\, \chi
        \label{eq:adjoint_equation_electrons_ions}
        \\
		&\lambda^\pm\vert_{t=T} = - \varphi^\pm
		\label{eq:adjoint_equation_initialcondition}
	\end{align}
\end{subequations}
with
\[
     \Phi_{f_B,\lambda}^\pm \coloneqq  \theta^\pm \pm \nabla_x\cdot (-\Delta_x)^{-1} \left(\sum_\pm \int_{\Omega_v} \mu_v^\pm \lambda^\pm \partial_v f_B^\pm\,\mathrm dv \right).
\]

As the adjoint system is linear, we proceed as in \cite[Section~5]{Knopf2018OCP_VP} to show that the system \eqref{eq:adjoint_model} possesses a unique strong solution $\lambda_B=(\lambda_B^-,\lambda_B^+)$ with 
\[ 
    \lambda_B^\pm \in H^1\big((0,T)\times\RR^D\big)\subset C\big([0,T];H^1(\RR^D)\big).
\]
Note that the cutoff function $\chi$ is necessary to obtain this well-posedness result. 
Since $\supp f^\pm(t) \subset \overline \Omega$ for all $t\in [0,T]$, we know that $\Phi_{f_B,\lambda}^\pm$ depends at most on the values of $\lambda_B^\pm$ in $[0,T]\times\overline\Omega$.
We can thus argue as in \cite[Remark~2(b)]{Knopf2018OCP_VP} to conclude that $\lambda_B^\pm\vert_{[0,T]\times\overline\Omega}$ is independent of the choice of the cutoff function $\chi$ provided that $R^*>0$ is chosen large enough.

In the subsequent theorem, we will see that also $\widetilde{J}'(B)$ depends at most on the values of $\lambda_B^\pm$ in $[0,T]\times\overline\Omega$. Consequently, $\widetilde{J}'(B)$ as well as the first-order necessary optimality condition are independent of the choice of $\chi$ if $R^*>0$ is chosen sufficiently large.

\begin{theorem} \label{OPTSYS}
    Let $B\in\Bad$ be arbitrary with corresponding state $f_{B}$,
    and let $\lambda_{B}^\pm = (\lambda_{B}^-,\lambda_{B}^+)$ be the corresponding strong solution of the adjoint system \eqref{eq:adjoint_model}.
    
    Then, the Fr\'echet derivative of the reduced cost functional $\widetilde{J}$ at the point $B$ is given by
    \begin{align}
	\widetilde{J}'(B)[H] = \alpha (B,H)_{\mathcal V}
	+  \big( \mathcal G(B) , H \big)_{L^2([0,T]\times \Omega_x)}
	\label{eq:reduced_Gradient}
    \end{align}
    for all directions $H\in \mathcal V$, where
    \begin{align}
         \mathcal G(B) \coloneqq  - \sum_\pm \mu_x^\pm\mu_v^\pm \int_{\Omega_v} v \otimes \partial_v f_{B}^\pm \lambda_{B}^\pm\,\mathrm dv.
         \label{eq:definition_rhs_gradient}
    \end{align}
    Now, suppose that $B^*\in \Bad$ is a locally optimal solution of the optimal control problem \eqref{eq:reduced_optim_problem}, i.e., there exists $\delta>0$ such that      
    \begin{align*}
        \widetilde{J}(B^*) \le \widetilde{J}(B) \quad\text{for all $B\in\Bad$ with $\|B-B^*\|_{\mathcal V} < \delta$}.
    \end{align*}
    Then, $B^*$ satisfies the variational inequality
    \begin{align}
        \label{eq:Optimality_condition}
        \widetilde{J}'(B^*)[B-B^*] 
        = \alpha (B^*,B-B^*)_{\mathcal V}
	       +  \big( \mathcal G(B^*) , B-B^* \big)_{L^2([0,T]\times \Omega_x)}
        \ge 0
    \end{align}
    for all $B\in\Bad$. In particular, if $a_{0j} < \partial_t B^*_j < b_{0j}$, $a_{ij} < \partial_{x_i} B^*_j < b_{ij}$ a.e.~in $[0,T]\times \Omega_x$ for all $i\in\{1,...,d_x\}$ and $j\in\{1,...,d_B\}$, $B^*$ is the weak solution of the following elliptic system with homogeneous Neumann boundary conditions:
    \begin{subequations}\label{eq:BVP}
	\begin{alignat}{2}
        \label{eq:BVP:1}
        - \alpha\,\mathrm{div}_{(t,x)} \big( \mathbf{K} \nabla_{(t,x)} B^*_j \big) + \alpha\,B^*_j &= - [\mathcal G(B^*)]_j
        &&\quad\text{in $(0,T)\times \Omega_x$},
        \\
		\label{eq:BVP:2}
        \big( \mathbf{K} \nabla_{(t,x)} B^*_j \big) \cdot \mathbf{n} &= 0
        &&\quad\text{on $\partial\big((0,T)\times \Omega_x\big)$}
	\end{alignat}
    for all $j\in\{1,...,d_B\}$.
    Here, $\mathbf{K} = \mathrm{diag}(\kappa_t,\kappa_x,...,\kappa_x) \in \RR^{(1+d_x)\times(1+d_x)}$ and $\mathbf n$ denotes the outer unit normal vector field on $\partial\big((0,T)\times \Omega_x\big)$.
    \end{subequations}
\end{theorem}

\medskip

\begin{proof}
Let now $B\in\Bad$ and $H\in\mathcal V$ be arbitrary.
Employing the chain rule, we obtain
\begin{align*}
    \widetilde{J}'(B)[H] 
    &= \alpha (B,H)_{\mathcal V}
	+ \sum_\pm \int_{\Omega_x} \int_{\Omega_v} \varphi^\pm\,  f_B^{\prime\pm}[H](T)\,\mathrm dv\mathrm dx
    \\ 
    &\quad +
    \sum_\pm  \int_0^T \int_{\Omega_x} \int_{\Omega_v} \theta^\pm\, f_B^{\prime\pm}[H]\,\mathrm dv\mathrm dx \mathrm dt.
\end{align*}
As $\lambda_B$ is the corresponding strong solution of the adjoint system, we can replace $\theta^\pm$ and $\varphi^\pm$ by means of \eqref{eq:adjoint_equation_electrons_ions} and \eqref{eq:adjoint_equation_initialcondition}, respectively. This yields
\begin{align*}
    &\widetilde{J}'(B)[H] 
    = \alpha (B,H)_{\mathcal V}
	- \sum_\pm \int_{\Omega_x} \int_{\Omega_v}  \lambda_B^\pm(T) \, f_B^{\prime\pm}[H](T)\,\mathrm dv\mathrm dx
    \\ \notag 
    &\quad +
    \sum_\pm \int_0^T \int_{\Omega_x} \int_{\Omega_v} \big[ \partial_t \lambda_B^\pm + \mu_x^\pm v\cdot\partial_x \lambda_B^\pm +\mu_v^\pm (E_{f_B}+\mu_x^\pm v\otimes B)\cdot\partial_v \lambda_B^\pm \big]\, f_B^{\prime\pm}[H]\,\mathrm dv\mathrm dx \mathrm dt
    \\ \notag 
    &\quad + \sum_\pm \int_0^T \int_{\Omega_x} \int_{\Omega_v} \mp \nabla_x\cdot (-\Delta_x)^{-1} \left(\sum_\pm\int_{\Omega_v} \mu_v^\pm \lambda_B^\pm \partial_v f_B^\pm\,\mathrm dv \right) \, f_B^{\prime\pm}[H]\,\mathrm dw\mathrm dx \mathrm dt .  
\end{align*}
We now recall that 
\[
    E_{f_B'[H]} = -\nabla_x (-\Delta_x)^{-1} \left( \sum_\pm \pm \int_{\Omega_v} f_B^{\prime\pm}[H] \, \mathrm dv\right)
\]
and that $f_B^{\prime\pm}[H]\vert_{t=0} = 0$.
Hence, by applying integration by parts, we obtain
\begin{align*}
    &\widetilde{J}'(B)[H] 
    = \alpha (B,H)_{\mathcal V}
    \\ \notag 
    &\quad 
    - \sum_\pm \int_0^T \int_{\Omega_x} \int_{\Omega_v} \big[ \partial_t f_B^{\prime\pm} + \mu_x^\pm v\cdot\partial_x f_B^{\prime\pm} +\mu_v^\pm (E_{f_B}+\mu_x^\pm v\otimes B)\cdot\partial_v f_B^{\prime\pm} \big] \lambda_B^\pm \,\mathrm dv\mathrm dx \mathrm dt
    \\ \notag 
    &\quad 
    - \sum_\pm \int_0^T \int_{\Omega_x} \int_{\Omega_v} \big[ \mu_v^\pm E_{f_B'} \cdot \partial_v f_B^\pm \big] \lambda_B^\pm \,\mathrm dv\mathrm dx \mathrm dt,  
\end{align*}
where we wrote $f_B^{\prime\pm}$ instead of $f_B^{\prime\pm}[H]$ for brevity.
Since $f_B'[H]$ is a strong solution of \eqref{VP:lin}, we can use \eqref{VP:lin} to reformulate the last two lines of the above equation. We thus have
\begin{align*}
    &\widetilde{J}'(B)[H] 
    = \alpha (B,H)_{\mathcal V} - \int_0^T \int_{\Omega_x} \left( \sum_\pm \mu_x^\pm\mu_v^\pm \int_{\Omega_v}v\otimes \partial_v f_B^\pm \lambda_B^\pm\,\mathrm dv \right) \cdot H \,\mathrm dx\mathrm dt.
\end{align*}
This verifies \eqref{eq:reduced_Gradient}.

Now suppose that $B^*$ is a \textit{locally optimal control} for the optimal control problem \eqref{eq:reduced_optim_problem}.
Since $\Bad$ is convex, we infer
\begin{align*}
	0\le \widetilde{J}'(B^*)[B-B^*] = \alpha (B^*,B-B^*)_{\mathcal V}
	       +  \big( \mathcal G(B^*) , B-B^* \big)_{L^2([0,T]\times \Omega_x)}
\end{align*}
for all $B\in\Bad$, where the equality holds due to \eqref{eq:reduced_Gradient}. This proves \eqref{eq:Optimality_condition}.

In particular, if $a_{0j} < \partial_t B^*_j < b_{0j}$, $a_{ij} < \partial_{x_i} B^*_j < b_{ij}$ a.e.~in $[0,T]\times \Omega_x$ for all $i\in\{1,...,d_x\}$ and $j\in\{1,...,d_B\}$, the function $B^*$ does not hit the box constraints and we can thus even vary in every direction. 
We thus obtain 
\begin{align*}
	0= \widetilde{J}'(B^*)[H] = \alpha (B^*,H)_{\mathcal V}
	       +  \big( \mathcal G(B^*) , H \big)_{L^2([0,T]\times \Omega_x)}
\end{align*}
for all directions $H\in\mathcal V$.
As this is exactly the weak formulation of \eqref{eq:BVP}, we conclude that $B^*$ is a weak solution of this elliptic system. Hence, the proof is complete.
\end{proof}

\section{Monte Carlo framework and numerical optimization}
\label{sec:MC_framework}

Our aim is to solve the optimality system 
consisting of the \textit{state equations} \eqref{VP:PM:d}, the \textit{adjoint equations} \eqref{eq:adjoint_model} and the \textit{first-order necessary optimality condition} \eqref{eq:Optimality_condition} by means of a
gradient-based optimization scheme, where the state equations and the adjoint equations are solved by a Particle-In-Cell method that is based on a Monte Carlo framework. 
Recall that the main purpose of the adjoint variable $\lambda^\pm$ is to allow for an efficient computation of the reduced gradient. 

To compute the adjoint variable, we also interpret $\lambda^\pm$ as a distribution function and we solve the adjoint equations by a generalized Particle-In-Cell method.
In view of our implementation, we assume that both $\theta^\pm$ and $\varphi^\pm$ can be interpreted as a probability distribution from which random numbers obeying this distribution can easily be sampled. We specify the following structure of $\theta^\pm$ and $\varphi^\pm$, that is consistent with \cref{asmpt:potentials}:
\begin{subequations}
	\label{def:phitheta}
	\begin{align}
		\theta^\pm(t,z) &= - \frac{C_{\theta^\pm}}{\sqrt{(2\pi)^{D}\det(\Sigma_{\theta^\pm})}} \exp\left(-\frac{1}{2}(z-z^\pm_\mathfrak{d}(t))^T\Sigma^{-1}_{\theta^\pm}(z-z^\pm_\mathfrak{d}(t)) \right), 
		\\
		\varphi^\pm(z) &= -\frac{C_{\varphi^\pm}}{\sqrt{(2\pi)^{D}\det(\Sigma_{\varphi^\pm})}}\exp \left( -\frac{1}{2}(z-z^\pm_T)^T\Sigma^{-1}_{\varphi^\pm}(z-z^\pm_T) \right), 
	\end{align}
\end{subequations}
with $C_{\theta^\pm}, C_{\varphi^\pm} >0$,
where $\Sigma_{\theta^\pm}$ and $\Sigma_{\varphi^\pm}$ are co-variance matrices in $\mathbb{R}^{D\times D}$ that we assume to be diagonal matrices. Moreover, the functions $z^\pm_\mathfrak{d}(t)$ and $z^\pm_T$ represent desired phase-space profiles for the numerical particles and need to be adjusted according to the optimization goal.

In this section, we focus on the case $d=(d_x,d_v)=(1,2)$ and we choose $\Omega_x$ to be a bounded interval.
However, we point out that our methodology also applies analogously in higher dimensions. 
To restate the optimality system in dimension $d=(1,2)$, we use the relations $v\cdot\partial_x=v_1\partial_x$, $E\cdot\partial_v=E\partial_{v_1}$, $v\otimes B=(v_2B,-v_1B)$, $(v\otimes B)\cdot\partial_v=v_2B\partial_{v_1}-v_1B\partial_{v_2}$, $-\nabla_x(-\Delta_x)^{-1}a=-\nabla_x\cdot(-\Delta_x)^{-1}(a,b)^T=\int_{-\infty}^xa\,\mathrm dy$ for suitable functions $a$ and $b$. 
We further assume that the functions $-a_{ij}$ and $b_{ij}$ in the definition of $\Bad$ are chosen sufficiently large such that all controls $B$ that occur in our optimization procedure remain strictly between the box constraints. Then, according to \cref{OPTSYS}, the first-order necessary optimality condition is actually an equation.
Thus, the optimality system in dimension $(1,2)$ consists of the state equations
\begin{subequations}\label{eq:forward_model_1D2V}
	\begin{align}
		&\partial_tf^\pm + \mu_x^\pm v_1\partial_xf^\pm +\mu_v^\pm(E_f+\mu_x^\pm v_2B)\partial_{v_1}f^\pm - \mu_x^\pm \mu_v^\pm v_1B\partial_{v_2}f^\pm=0,\\
        &f^\pm \vert_{t=0} = \mathring{f}^\pm
		\label{eq:poisson_1D2V}
	\end{align}
\end{subequations}
with 
\begin{align}
    \label{eq:electric_field_1D2V}
    E_f=\int_{-\infty}^x\rho_f\,dy=\int_{-\infty}^x\int_{\Omega_v}(f^+-f^-)\,\mathrm dv\mathrm dy,
\end{align}
the adjoint equations
\begin{subequations}\label{eq:adjoint_model_1D2V}
	\begin{align}
	&\begin{aligned}
		\partial_t\lambda^\pm + \mu_x^\pm v_1\partial_x \lambda^\pm +\mu_v^\pm (E_f + \mu_x^\pm v_2 B)\partial_{v_1}
		\lambda^\pm - \mu_x^\pm \mu_v^\pm v_1B \partial_{v_2}\lambda^\pm
		\\ = \theta^\pm\mp \int_{-\infty}^x \left(\sum_\pm \int_{\Omega_v} \mu_v^\pm \lambda^\pm \partial_{v_1} f^\pm\,\mathrm dv \right)\mathrm dy,
	\end{aligned}\\
	&\lambda^\pm \vert_{t=T}=-\varphi^\pm,
	\end{align}
\end{subequations}
and the first-order necessary optimality condition
\begin{subequations}\label{eq:BVP*}
	\begin{alignat}{2}
        \label{eq:BVP:1*}
        - \alpha\,\mathrm{div}_{(t,x)} \big( \mathbf{K} \nabla_{(t,x)} B \big) + \alpha\,B &= - \mathcal G(B)
        &&\quad\text{in $(0,T)\times \Omega_x$},
        \\
		\label{eq:BVP:2*}
        \big( \mathbf{K} \nabla_{(t,x)} B \big) \cdot \mathbf{n} &= 0
        &&\quad\text{on $\partial\big((0,T)\times \Omega_x\big)$}
	\end{alignat}
 \end{subequations}
for all $j\in\{1,...,d_B\}$,
where $\mathcal G(B)$, $\mathbf{K}$ and $\mathbf{n}$ are defined as in \cref{OPTSYS}. Recall that we now have $d_x = 1$ and $d_B = 1$ meaning that $B$ and $\mathcal G(B)$ are scalar functions.

We now proceed with the presentation of our numerical scheme. 
In general, the time evolution of gases can be described by means of particle methods (see, e.g., \cite{Crowe2011multiphaseflowsParticles}). In particular, such methods are also applicable to describe plasmas.
In this sense, a Monte Carlo framework can be used to exploit the inherent natural structure of gases or plasmas.
However, there are also other possibilities to solve the equations appearing in this paper. Some of them are described, for instance, in \cite{Tajima}.

In the present paper, we implement a so-called Particle-In-Cell (PIC) method (see, e.g., \cite{Birdsall2004PlasmaPhysicsSimulation}).
The idea is to approximate the distribution functions by numerical particles which are represented by labelled pointers to structures that contain the physical information such as position and velocity.
The particles are considered in the Lagrangian framework, i.e., the evolution of the particles is computed by solving the characteristic systems of the Vlasov equations. We point out that this would not require the introduction of a mesh to discretize the considered region in phase space.
However, since electric and magnetic fields are involved in the characteristic system, we further need to compute these quantities and for this purpose, a mesh actually needs to be generated.
The fields are then computed at the grid points of the mesh. For any numerical particle, the influence of the fields on this particle is calculated via interpolation techniques.

Essentially, a timestep in the PIC procedure consists of changing the content of the structure representing a particle,
namely the position and velocity, and adding or subtracting pointers (numerical particles), if required. 
More precisely, we need to solve the characteristic systems 
\begin{align}
    \label{eq:ODE_system_Boris}
    \dot{x}(t) = \mu^\pm_x v_1(t) ,
    \quad
    \dot{v_1}(t) &= \mu^\pm_v\, E_f(t,x(t)) + \mu^\pm_x\mu^\pm_v v_2(t)  \,B(t,x(t)), 
    \\ \notag
    \dot{v_2}(t) &= - \mu^\pm_x\mu^\pm_v v_1(t)  \,B(t,x(t))		
\end{align}
by a suitable ODE solver. 
In this context, also the electric field $E_f$ needs to be computed by a numerical method.

The free streaming time is given by $\Delta t$, which is chosen as an input parameter.
Within this time-lapse, the microscopic equations of motion \eqref{eq:ODE_system_Boris} have to be integrated. For this purpose, we apply the so-called \textit{Boris pusher} (see, e.g.,  \cite{Boris1970cylrad, Zenitani2018}).
An important property of this symplectic method is that it preserves the microscopic energy of the numerical particles even in the presence of a magnetic field.
To formulate the algorithm, we extend the position and velocity to vectors in $\RR^3$ by filling them up with zeros. This has the advantage that the algorithm can easily be adapted to higher dimensional situations.
The implementation of the Boris pusher is then given by \autoref{algo:Boris_Pusher}. 
\begin{algorithm}
	\caption{Boris pusher ($d=(1,2)$)}
	\label{algo:Boris_Pusher}
	\begin{algorithmic}[1]
		\REQUIRE velocity and position $v, x \in \RR^3$, magnetic field $B \in \RR$, electric field $E \in \RR$, charge $q$, timestep $\Delta t$
		\ENSURE updated position and velocity
		\STATE Define $\boldsymbol{E} = (E,0,0) \in \RR^3$
		and $\boldsymbol{B} = (0,0,B) \in \RR^3$
		\STATE Calculate $r = q \boldsymbol{B} \frac{\Delta t}{2}$
		and $s = 2r/(1+||r||_2^2)$
		\STATE Calculate $v^- = v + q\boldsymbol{E}\frac{\Delta t}{2}$, then $v' = v^- + v^- \times r$ and lastly $v^+ = v^- + v' \times s$
		\STATE Update the velocity $v = v^+ + q\boldsymbol{E} \frac{\Delta t}{2}$ and the position $x = x + v \Delta t$
		\RETURN $(x,v)$
	\end{algorithmic}
\end{algorithm}

While updating the position, one has to take the boundedness of the physical domain and the boundary condition into account.
Throughout the following part of our numerical work, we consider periodic boundary conditions with respect to the position coordinate. 

To initialize the numerical particles whose movement is then determined by \autoref{algo:Boris_Pusher}, we use a Monte Carlo approach.
This means we discretize the initial data $\mathring{f}^\pm$ by initializing 
a large number $N_{f}$ of numerical particles.
Moreover, we consider a partition of the time interval $[0, T]$ into $N_t$ subintervals of size $\Delta t = T/N_t$. 
The time of the $k$-th time step is thus given by $t^k = k \Delta t$, $k=0,\ldots,N_t$.

In our implementation, we define $F$ as the list of labelled pointers to structures that represent the numerical particles. For any $p=1,\ldots,N_f$,
we denote by $F_\pm^k[p]$ the pointer to the $p$-th particle of the ions or electrons, respectively, at the $k$-th timestep.
Further, let $F_\pm^k[p].v_i$, $i=1,2$, be the velocity of the $p$-th particle at the $k$-th timestep, and let $F_\pm^k[p].x$ be the position of the  $p$-th particle at the $k$-th timestep. 
Analogously, $\Lambda$ denotes the list of labelled pointers to structures representing the adjoint particles.

To initialize a list of particles using a density from which it is possible to directly generate samples (e.g., Gaussian or uniform density), we sample the position and velocity of $N_f$ particles according to the given distribution.

However, for several reasons that will be explained below, it may be desirable to have an initial condition given by some arbitrary distribution $g=g(x,v)$.
To realize this, there are several possibilities.
For instance, one can use the so-called Inverse Transform Sampling technique for which the inverse of the cumulative distribution function has to be known (see, e.g., \cite[Chapter 2.1]{Steele1987nonUniform}).
In this paper, we employ an acceptance-rejection technique to sample from a desired probability density $g=g(x,v)$ with the help of another probability density $h=h(x,v)$ (cf.~\cite[Algorithm $A_1$]{Casella2004AcceptanceRejection}). 
In this case, $h$ is to be chosen in such a way that drawing samples from $h$ is straightforward to implement.
Moreover, we demand that there exists a constant $k>0$ such that $g(x,v)\leq k \, h(x,v)$ for all $(x,v) \in \Omega$.
The pseudo-code of this method is presented in \autoref{algo:initialize_arbitrary}.
\begin{algorithm}
	\caption{Acceptance-Rejection-Algorithm to approximate general initial data}
	\label{algo:initialize_arbitrary}
	\begin{algorithmic}[1]
		\REQUIRE desired density $g(x,v)$, helper density $h(x,v)$, desired number of particles $N$
		\STATE Set $p=1$; Initialize empty list $L^0$
		\WHILE{$p < N$}
			\STATE Generate sample $y_p$ from $h$ 
			\STATE Generate sample $u_p$ from uniform distribution $\mathcal{U}(0,1)$
			\IF{$u_p < \frac{g (y_p)}{k\, h(y_p)} $ }
				\STATE accept sample $y_p$ and generate particle with corresponding position and velocity and add it to the list $L^0$
				\STATE $p = p+1$
			\ENDIF \quad // else: reject sample
		\ENDWHILE
		\RETURN $L^0$
	\end{algorithmic}
\end{algorithm}
As we need to compute the electric field, the right-hand sides of the adjoint equations as well as the gradient of our cost functional during the optimization process, we need to introduce a mesh that discretizes the considered region in phase space.

For this purpose, we consider a bounded domain of velocities $\Omega_v \coloneqq (-v_{\max},v_{\max})^2 \subset \RR^2$. Here, $v_{\max}>0$ serves as a threshold that the velocity components of the numerical particles shall not exceed.
The setting of this parameter makes sense whenever initial distributions $\mathring{f}^\pm=\mathring{f}^\pm(x,v)$ that are either compactly supported or decay to zero at a sufficiently high rate as $|v|\to\infty$ are considered. 
In particular, this allows for initial distributions that are Gaussian-like with respect to the velocity component.
Now, we define a partition of $\Omega_v$ in equally sized, non-overlapping, square cells of edge length $\Delta v = {2v_{\max}}/{N_v}$ with $N_v \geq 2$. 
On this partition, we consider a cell-centered representation of the velocities within
\begin{align*}
	\Omega_{\Delta v} \coloneqq \left\lbrace \; (v^l_1,v^m_2) \in \Omega_v \;\big\vert\; l,m = 1,\ldots,N_v \; \right\rbrace,
\end{align*}
where $v^l_1 \coloneqq \left(l-1/2\right)\Delta v - v_{\max}$ and $v^m_2 \coloneqq \left(m-1/2\right)\Delta v - v_{\max}$.

Similarly, we consider $\Omega_x = (0,x_{\max}) \subset \RR$ with $x_{\max}>0$, and we define a partition in intervals of length $\Delta x = x_{\max}/N_x$, where $N_x \geq 2$. 
On this partition, we define
\begin{align*}
	\Omega_{\Delta x} \coloneqq \left\lbrace \; x^i \in \Omega \;\big\vert\;  i = 1,\ldots,N_x \; \right\rbrace, \qquad x^i \coloneqq (i-1/2)\Delta x .
\end{align*}
We thus have the discretized phase-space $\Omega_{\Delta x} \times \Omega_{\Delta v}$.

Now, we denote by $(f^\pm)_{ilm}^k$ the occupation number of the cell with center $(x^i,v^l_1,v^m_2) \in \Omega_{\Delta x} \times \Omega_{\Delta v}$ at the $k$-th timestep. 
To construct this function, we count the particles at timestep $k$ that have position and velocity in the cell centered at $(x^i,v_1^l,v_2^m)$. 
Thus, we define 
\begin{align}
	(f^\pm)_{ilm}^k \coloneqq \sum_{p=1}^{N_f}  \mathbbm{1}_{ilm}
	\left( 
	F_\pm^k[p].x,F_\pm^k[p].v_1 ,F_\pm^k[p].v_2 \right),
	\label{eq:assembling_Distribution_f}
\end{align}
where $\mathbbm{1}_{ilm}(\cdot,\cdot,\cdot)$ denotes the indicator function of the cell with center $(x^i,v_1^l,v_2^m)$. If a particle leaves the cell centered at $(x^i,v_1^l,v_2^m)$ and enters the cell centered at $(x^{i'},v^{l'}_1,v^{m'}_2)$, the value of $(f^\pm)_{ilm}$ is reduced by 1 and the value of $(f^\pm)_{i'l'm'}$ is increased by 1.
Notice that by choosing $v_{\max}$ large enough, the probability that the velocity of a particle exceeds the boundary of $\Omega_{v}$ is very low but possibly not zero. 
Therefore, in our implementation, we keep track of the velocity of the particles that exceed the boundary of $\Omega_{v}$ and provide a warning if the number of such particles surpasses a certain threshold.

To solve the adjoint system, we have to deal with velocities exceeding the threshold $v_{\max}$ in a different way.
As we need to add additional particles to account for the source terms (right-hand sides) in the adjoint kinetic equations, the phenomenon that particles will actually reach the bound $v_{\max}$ may actually occur very frequently.
However, due to \eqref{eq:reduced_Gradient}, the gradient of the cost functional depends on 
the values of the adjoint variables $\lambda^\pm(t)$ only on the support of $f^\pm(t)$. 
Therefore, as long as the support of $f^\pm(t)$ remains confined, we only need to consider adjoint particles with velocity components in the interior of $\Omega_v$ to compute the gradient. 
We compute the occupation numbers associated with the adjoint particles as 
\begin{align}
	(\lambda^\pm)_{ilm}^k = \sum_{p=1}^{N_{\lambda,k}^\pm}  \mathbbm{1}_{ilm}\left( \Lambda_\pm^k[p].x, \Lambda_\pm^k[p].v_1, \Lambda_\pm^k[p].v_2 \right),
	\label{eq:assembling_Distribution_q}
\end{align}
where $N_{\lambda,k}^\pm$ stands for the number of numerical particles at the $k$-th time step. 

To compute the electric field in dimension $(1,2)$, we choose the procedure given in \autoref{algo:Poisson_1D2V} that is a simple quadrature rule for integrating the density of the difference of the electron and ion density.
Therefore, we represent the electric field as
\begin{align}
	E_f(t,x) = \frac{1}{2}\left( \int_{-\infty}^x \rho_f(t,y)\, \mathrm d y - \int_{x}^{\infty}\rho_f(t,y)\,\mathrm d y	\right)
	\label{eq:Formula_ElectricField}
\end{align}
to circumvent high errors due to asymmetries made by our first-order approximation of the integrals.
This formula is justified by the identity
\begin{align*}
	\int_{-\infty}^x\rho_f(t,y)\,\mathrm d y 
    = \int_{-\infty}^\infty \rho_f(t,y)\,\mathrm d y - \int_{x}^{\infty}\rho_f(t,y)\,\mathrm d y
    = - \int_{x}^{\infty}\rho_f(t,y)\,\mathrm d y
\end{align*}
that is a consequence of quasi-neutrality $\int_{-\infty}^\infty \rho_f(t,y)\,\mathrm d y = 0$, which holds since we use the same amount of ions and electrons.
As we have several nested for-loops, this may need require a long computation time.
However, it is possible to parallelize the loop over $i$.
\begin{algorithm}
	\caption{Computation of the electric field (for $d=(1,2)$)}
	\label{algo:Poisson_1D2V}
	\begin{algorithmic}[1]
		\REQUIRE $((f^-)^k,(f^+)^k)$
		\ENSURE Value $E_f$ of integral in \eqref{eq:electric_field_1D2V}
		\FOR{$i=0$ \TO $N_x$} 
    		\FOR{$l,m = 0$ \TO $N_v$}
    		      \STATE Calculate charge density $\rho^k_{i} = \rho^k_{i} + ((f^+)^k_{ilm} - (f^-)^k_{ilm})$ 
    		\ENDFOR
    		\STATE Calculate $(E^+_f)^k_i = (E_f)^k_{i-1} + \rho^k_{i}\,\Delta x$ and $(E^-_f)^k_i = (E_f)^k_{N_x - i} + \rho^k_{N_x-i-1}\,\Delta x$ 
    		\STATE Assemble the electric field $(E_f)^k_i = \frac{1}{2}\left((E^+_f)^k_i - (E^-_f)^k_i\right)$ (cf. \eqref{eq:Formula_ElectricField})
		\ENDFOR
		\RETURN $E_f$
	\end{algorithmic}
\end{algorithm}

We can now present \autoref{algo:forward_solver} in order to solve the state equation \eqref{eq:forward_model_1D2V} using our Monte Carlo framework.
\begin{algorithm}
	\caption{Forward solver (for $d=(1,2)$)}
	\label{algo:forward_solver}
	\begin{algorithmic}[1]
		\REQUIRE Initial conditions $\mathring{f}^\pm$, control $B^k$ for $k=0,\ldots,N_t$
		\ENSURE time evolution of $f$
		\STATE Generate initial particles $F_\pm^0$ applying \autoref{algo:initialize_arbitrary} with $g=\mathring{f}^-$ and $g=\mathring{f}^+$ using the Gaussian density as helper density
		\FOR {$k=0$ \TO $N_t-1$}
			\STATE Calculate $E_f^k$ by \autoref{algo:Poisson_1D2V} with input $F_\pm^k$
			\FOR{$p=1$ \TO $N_f$}
				\STATE Update $F_\pm^k[p].x$ and $F_\pm^k[p].v$ applying the Boris pusher (\cref{algo:Boris_Pusher}) by means of linear interpolation of the electric field $E_f^k$ and the magnetic field $B^k$ 
			\ENDFOR
		\ENDFOR
	\end{algorithmic}
\end{algorithm}

\begin{rem}
    Notice that there are different procedures in the literature to solve the Vlasov--Poisson problem; see, e.g. \cite{Chertock2017GuideParticleMethods} for an overview. 
	In particular, there are schemes that directly use weighted particles and therefore avoid costly sampling procedures during their execution.
    Usually, when only the forward problem is to be solved numerically, the charge density $\rho_f$ is the actually important quantity, also because it can be visualized more easily than the distribution functions $f^\pm$. In this case, there are more efficient methods (e.g., projection methods) to directly recover the charge density from the numerical particles. In our case, however, we really need to compute the distribution functions $f$ as well as the adjoint variables $g$ on the same grid in order to use them for the computation of our control $B$.
    Therefore, to compute $f$ and $g$, we perform a piecewise constant interpolation of the corresponding numerical particles. Higher order interpolation (e.g., by splines) or scattered data approximation techniques would, of course, be more precise, but very expensive due to the high amount of numerical particles. Such techniques might be reasonable if only the forward problem is to be solved more accurately. However, in our situation, where we have to solve the forward and backward (adjoint) problem several times, it does not seem realizable given the current computation power. 
\end{rem}
Next, we present a similar scheme for solving the adjoint system. 
As we will see, we need to create additional numerical particles to account for the source terms on the right-hand sides of the adjoint Vlasov equations. 
The adjoint equations \eqref{eq:adjoint_model} can be rewritten as
\begin{align}
	\partial_t \lambda^\pm = \mathscr{F}^\pm(f^\pm,\lambda^\pm,B) + \mathscr{R}^\pm(\lambda^\pm,f^\pm) + \mathscr{S}^\pm.
	\label{eq:adjoint_structure}
\end{align}
Here, $\mathscr{F}^\pm$ represents the ``Vlasov'' terms that will be handled by the Boris pusher (\cref{algo:Boris_Pusher}).
The source term $\mathscr{S}^\pm$ is given by $\mathscr{S}^\pm = \theta^\pm$.
Moreover, $\mathscr{R}^\pm$ denotes the reaction term 
\begin{align}
	\mathscr{R}^\pm(\lambda^\pm,f^\pm) \coloneqq \mp \int_{-\infty}^x \int_{\Omega_v} \sum_\pm \mu_v^\pm \lambda^\pm \partial_{v_1} f^\pm \,\mathrm dv\mathrm dy.
    \label{eq:Reaction_term_definition}
\end{align}
Note that in \eqref{eq:Reaction_term_definition}, the derivative acts on the quantity $f$, which is already known as the forward problem is solved in advance. The derivative is approximated using a second-order finite difference scheme. We further notice that the term $\mathscr{R}^\pm$ is of the same type as the electric field.
Therefore, it can be solved by using the same integration method as for the computation of $E_f$.
As the reaction terms $\mathscr{R}^\pm$ appear on the right-hand side of the adjoint equations, there is no conservation of the total charge. For that reason, we might have to add additional numerical particles while solving the adjoint system.
We apply Euler's method along with an operator splitting to obtain the recursive formula
\begin{align}
	(\lambda^\pm)^{k+1}_{ilm} = (\lambda^\pm)^k_{ilm} 
    + \Delta t \, \Big( \mathscr{R}((\lambda^\pm)^k_{ilm},(f^\pm)^k_{ilm}) + \mathscr{S}^\pm \Big).
	\label{eq:Eulers_method_integral}
\end{align}
Hence, the number of particles in every cell centered at  $(x^i,v_1^l,v_2^m)$ in phase-space at timestep $k+1$ has to be changed according to the value of the second summand on the right-hand side of \eqref{eq:Eulers_method_integral}.
That is, particles have to be added (or possibly removed) in the phase-space cell centered at $(x^i,v_1^l,v_2^m)$.
The addition of particles is adding elements to the particle list with a given position and velocity. 
The removal of particles is more delicate since one needs to find a particle with position and velocity in the phase-space cell under consideration. This is of course very time-consuming. 
However, by choosing the constants $C_{\theta^\pm}$ and $C_{\varphi^\pm}$ in the definition of  $\theta^\pm$ and $\varphi^\pm$ (see \eqref{def:phitheta}) sufficiently large, we ensure that the right-hand sides of the adjoint equations remain non-negative. Then 
the removal of numerical particles is not necessary.
In \autoref{algo:linReact_creation}, we summarized the procedure in pseudo-code.
\begin{algorithm}
	\caption{Implementation of the reaction and source term at time $t^k$}
	\label{algo:linReact_creation}
	\begin{algorithmic}[1]
		\REQUIRE $\Lambda_\pm^k$, $F_\pm^k$, $\theta^\pm(t^k,z)$
		\ENSURE $\Lambda_\pm^k$ at the end of the timestep
		\STATE Calculate $(N^\pm)_\mathscr{R}^i$ as the value of the integral \eqref{eq:Reaction_term_definition} using Poisson solver Algorithm \ref{algo:Poisson_1D2V} solver for all cells $x^i$, $i=1,\ldots,N_x$.
		\FOR{ $i=1$ \TO $N_x$  }		
			\FOR{$l,m=1$ \TO $N_v$}
			\STATE Calculate the amount $N_{\theta^\pm}^{ilm}$ of particles that will be created in the phase-space cell $(x^i,v^l_1,v^m_2)$ using the analytic expression of $\theta^\pm(t,z)$ 
			\STATE Set $N_\pm^{ilm} = N_{\theta^\pm}^{ilm} - (N^\pm)_\mathscr{R}^i$
				\IF{$N_\pm^{ilm} > 0$}
					\STATE Create $N_\pm^{ilm}$ particles with uniform distributed velocity in the current velocity-cell $(v_1^l,v_2^m)$ and uniform distributed position in current position-cell $x^i$
					\STATE Add them to $\Lambda_\pm^k$
				\ENDIF // if $N_\pm^{ilm} \leq 0$ do nothing
			\ENDFOR
		\ENDFOR
		\STATE Add generated particles to the existing ones in $\Lambda_\pm^k$
	\end{algorithmic}
\end{algorithm}
Gathering all subroutines, we can form the final method to solve the backward model, see \autoref{algo:backward_solver}.
Inside this algorithm, we also take the influence of the source term $\mathscr{S}^\pm$ into account.
Note that the approach in \autoref{algo:backward_solver} is of first-order since we apply a first-order operator splitting. 
\begin{algorithm}
	\caption{Backward solver (for $d=(1,2)$)}
	\label{algo:backward_solver}
	\begin{algorithmic}[1]
		\REQUIRE Functions $\varphi^\pm$ and $\theta^\pm$, control $B^k$,  $k=0,\ldots,N_t$, solution and resulting electric field $(F_\pm^k,E^k_f)$, $k=0,\ldots,N_t$ of the forward problem
		\ENSURE time evolution of $\lambda$
		\STATE Generate particles for the terminal condition 
		$(\Lambda^{N_t}_\pm)$ using $-\varphi^\pm$ and direct sampling
		\FOR {$k=N_t-1$ \TO $0$}
		\FOR{$p=1$ \TO $N_q^k$}
		\STATE Update $\Lambda_\pm^k[p].x$ and $\Lambda_\pm^k[p].v$ according to Boris algorithm
		\ENDFOR
		\STATE Solve the reaction and source term using $F_\pm^k$, $\theta^\pm$ and \autoref{algo:linReact_creation}.
		\ENDFOR
	\end{algorithmic}
\end{algorithm}

Next, we focus on the gradient equation \eqref{eq:reduced_Gradient}. 
For the update of the control, we need to calculate the representative of the gradient of the reduced cost functional in the Hilbert space $\mathcal V$. 
This representation can be derived as follows. 
Assuming enough regularity of the reduced cost functional, the first-order Taylor expansion of $\widetilde{J}(B)$ in a Hilbert space $X$ reads as
\begin{align*}
	\widetilde{J} (B+\epsilon \, \delta B) =  \widetilde{J}(B) +
	\epsilon \, \left(\nabla \widetilde{J}(B) |_X, \delta B \right)_X 
	+ \BigO{\epsilon^2}.
\end{align*} 
Here, the gradient $\nabla \widetilde{J}(B) |_X$ depends on the choice of the Hilbert space $X$ and its corresponding inner product. Choosing $X=L^2((0,T)\times\Omega_x)$ and proceeding as in the derivation of \eqref{eq:BVP}, we obtain
\begin{align}
    \label{eq:DJL2x}
	\nabla \widetilde{J}(B) |_{L^2}
    &= - \alpha\,\mathrm{div}_{(t,x)} \big( \mathbf{K}  \nabla_{(t,x)} B \big) 
    + \alpha B + \mathcal G(B)
    \\ \notag
    &= -\alpha \kappa_t \partial_t^2 B - \alpha \kappa_x \partial_x^2 B + \alpha B
    + \mathcal G(B)
\end{align}
For the choice $X=\mathcal{V}$, we can determine the gradient based on the fact that the Taylor series must be identical term-by-term regardless of the choice of $X$, i.e.,
\begin{align*}
\left( \nabla \widetilde{J}(B) |_{\mathcal V}, \delta B \right)_{\mathcal{V}}=
\left( \nabla \widetilde{J}(B) |_{L^2}, \delta B \right)_{L^2} .
\end{align*}
Recalling the definition of the inner product on $\mathcal V$ (see \eqref{def:Vscp}), we proceed as in the derivation of \eqref{eq:BVP} to deduce that $\nabla \widetilde{J}(B) |_{\mathcal V}$ is given by the elliptic boundary value problem
\begin{subequations}
 \label{eq:gradH1x} 
\begin{align}
	- \kappa_t\, \partial_t^2 \nabla \widetilde{J}(B) |_{\mathcal V}  
	- \kappa_x\, \partial_x^2 \nabla \widetilde{J}(B) |_{\mathcal V} 
	+ \nabla \widetilde{J}(B) |_{\mathcal V}
	&=
	\nabla \widetilde{J}(B) |_{L^2}
    &&\quad\text{in $(0,T)\times \Omega_x$},
    \\
    \big(\kappa_t\, \partial_t \nabla \widetilde{J}(B) |_{\mathcal V} , 
        \kappa_x\, \partial_x \nabla \widetilde{J}(B) |_{\mathcal V}\big)^T
    \cdot \mathbf{n}  &= 0
    &&\quad\text{on $\partial\big((0,T)\times \Omega_x\big)$},
\end{align}
\end{subequations}
where $\mathbf n$ stands for the outer unit normal vector field on $\partial\big((0,T)\times \Omega_x\big)$.

To compute the gradient $\nabla \widetilde{J}(B) |_{L^2}$ according to \eqref{eq:DJL2x}, we employ a finite-difference discretization.
We assemble the optimization gradient in a vector $g \in \RR^{(N_x-2)}$. 
Therefore, a rectangular quadrature rule is used to approximate the integrals in the definition of $\mathcal G(B)$ (see \eqref{eq:definition_rhs_gradient}). In this way, we obtain for $k = 0,\ldots,N_t$ and $i=2,\ldots,N_x-1$ the formula
\small
\begin{align*}
	G_i^k \coloneqq -  \, 
	\sum_\pm \mu^\pm_x\mu^\pm_v\sum_{l,m}^{N_v-1} \left(
	l (f^\pm)_{i,l,m+1}^k
    -(l-m)(f^\pm)_{ilm}^k
	-
	m (f^\pm)_{i,l+1,m}^k
	  \right) (\lambda^\pm)_{ilm}^k .
\end{align*}
\normalsize
In view of \eqref{eq:gradH1x}, the finite-difference approximation of the gradient $\nabla \widetilde{J}(B) |_{L^2}$ reads as
 \small
\begin{align*}
	g^k_i \coloneqq \alpha B^k_i 
	- \frac{\alpha\,\kappa_x}{\Delta x^2}\left( B^k_{i+1}-2B^k_i+B^k_{i-1} \right) 
	- \frac{\alpha\,\kappa_t}{\Delta t^2}\left( B^{k+1}_i-2B^k_i+B^{k-1}_i \right) 
	+  (\Delta v)^2(\Delta t) G^k_i
\end{align*}
\normalsize
for $ k = 0,\ldots,N_t$ and $i=2,\ldots,N_x-1$.

Now, the gradient $\nabla \widetilde{J}(B) |_{\mathcal V}$ is computed by solving the elliptic problem \eqref{eq:gradH1x} by means of a standard finite-difference method. 
The resulting linear system is then given by a pentadiagonal block matrix and can be solved efficiently by a generalized Thomas algorithm.

Keep in mind that we needed to compute the gradient in the correct Hilbert space $\mathcal{V}$ in order to be consistent with the formulation of our optimal control problem \eqref{eq:reduced_optim_problem}.
Furthermore, recall from Lemma~\ref{lem:bad} that the choice of $\mathcal{V}$ leads to Lipschitz continuous controls.

To find a locally optimal control, we combine deterministic gradient-based optimization techniques with our particle methods.
With the preparation above, we can now formulate the algorithm that provides the gradient $\nabla \widetilde{J}(B) |_{\mathcal V}$ which is required in our optimization scheme. 
\begin{algorithm}
	\caption{Calculate the gradient $\nabla \widetilde{J}(B) |_{\mathcal V}$}
	\label{algo:CalculateGradient}
	\begin{algorithmic}[1]
		\REQUIRE control $B(t,x)$, $\mathring{f}^\pm(x,v)$, $\varphi^\pm(x,v)$, $\theta^\pm(t,x,v)$
		\STATE Solve the Vlasov--Poisson system \eqref{eq:forward_model_1D2V} using Algorithm \ref{algo:forward_solver} 
		with inputs $\mathring{f}^\pm$, $B$
		\STATE Solve adjoint system \eqref{eq:adjoint_model_1D2V} using Algorithm \ref{algo:backward_solver} 
		with inputs $-\varphi^\pm$, $\theta^\pm$, $B$
		\STATE Compute the distribution functions $f^\pm$ and $\lambda^\pm$ according to \eqref{eq:assembling_Distribution_f} 
		and \eqref{eq:assembling_Distribution_q}
		\STATE Assemble $\nabla \widetilde{J}(B) |_{L^2}$ according to \eqref{eq:DJL2x}
		\STATE Compute $\nabla \widetilde{J}(B) |_{\mathcal V}$ solving \eqref{eq:gradH1x}
	\end{algorithmic}
\end{algorithm}
We remark that, based on \autoref{algo:CalculateGradient}, we can implement many different gradient-based optimization schemes to solve our optimization problem (see, e.g., \cite{BorziSchulz2011}). 
In our case, we choose the non-linear conjugate gradient (NCG) method. For extensions of this method to the stochastic setting, we refer to \cite{Schraudolph2002towardsStochasticNCG, JinZhang2019StochasticCG}. 
This is an iterative method that constructs a minimizing sequence of 
control functions $(B^\ell)$ by the following algorithm.
\begin{algorithm}
	\caption{NCG scheme}
	\label{algo:OptAlgo}
	\begin{algorithmic}[1]
		\REQUIRE $B^0(t,x)$, $\mathring{f}^\pm(x,v)$, $\varphi^\pm(x,v)$, $\theta^\pm(t,x,v)$
		\STATE $\ell = 0$, $E > tol$
		\STATE Compute $h^0 = -\nabla \widetilde{J}(B^0)|_{\mathcal V}$ using Algorithm \ref{algo:CalculateGradient}
		\WHILE{$E >tol$ \AND $\ell<\ell_{\max}$}
		\STATE Use a line-search scheme to determine the step-size $\sigma_\ell$ along $h^\ell$
		\STATE Update control: $B^{\ell+1} = B^\ell + \sigma_\ell \, h^\ell$
		\STATE Compute $d^{\ell+1}= \nabla \widetilde{J} (B^{\ell+1}) |_{\mathcal V}$ using Algorithm \ref{algo:CalculateGradient}
		\STATE Compute $\beta_\ell$ using the Fletcher--Reeves formula
		\STATE Set $h^{\ell+1} = - d^{\ell+1} + \beta_\ell \, h^\ell$
		\STATE $E=\|B^{\ell+1}-B^\ell\|_2$
		\STATE Set $\ell = \ell+1$
		\ENDWHILE
		\RETURN $(B^\ell,f^\ell)$
	\end{algorithmic}
\end{algorithm}
In \autoref{algo:OptAlgo}, the tolerance $tol>0$ and the maximum number of iterations $\ell_{\max} \in \NN$ are used as termination criteria. 
We use backtracking line-search with Armijo condition. 
The factor $\beta_\ell$ is based on the Fletcher--Reeves formula; see \cite{BorziSchulz2011} for more details and references.

An estimate of the computational cost of one iteration of \cref{algo:OptAlgo} can be calculated as follows.
The cost for solving the Vlasov--Poisson equation using \cref{algo:forward_solver} is the sum of $\BigO{N_f}$ for generating the initial condition, and $\BigO{N_t N_x N_v^2 N_f^2}$ since in every timestep, we have to assemble the probability density functions ($\BigO{N_f}$), calculate the electric field according to \cref{algo:Poisson_1D2V} ($\BigO{N_xN_v^2}$) and calculate for every particle a new position and velocity ($\BigO{N_f}$). 
Hence we have an estimated cost of $\BigO{N_tN_xN_v^2N_f^2}$ for solving the Vlasov--Poisson model equation.
For solving the adjoint equation using \cref{algo:backward_solver}, we do not need to solve the equation for the electric field but need to include the source- and reaction-term with \cref{algo:linReact_creation} for which we estimate the cost to be $\BigO{N_x N_v^2}$. 
This leads to an estimated cost of $\BigO{N_t(N_f^2 + N_xN_v^2)}$ for solving the adjoint model.
Furthermore, we have $\BigO{N_tN_x}$ operations for the integration required for calculating the $L^2$ gradient in time and space and $\BigO{N_tN_x}$ to compute the $H^1$ gradient using the Thomas algorithm.
In summary, the computational complexity of a single optimization iteration is $\BigO{N_tN_xN_v^2N_f^2 + N_tN_x}$.
This estimate is without considering the linesearch, for which $\BigO{N_t N_x N_v^2}$ operations are required for calculating the functional and $\BigO{N_t N_x N_v^2 N_f^2}$ for solving the Vlasov--Poisson equation for this purpose.

\begin{rem}
    We point out that to be consistent with the mathematical analysis in Section~\ref{sec:Analysis_OCP}, one would actually have to include a projection onto the set $\mathbb{B}_{\mathrm{ad}}$ where bounds on the derivative of the controls are prescribed.
    However, if $-a_{ij}$ and $b_{ij}$ in the definition of $\mathbb{B}_{\mathrm{ad}}$ (see \cref{def:admissible_controls}) are a priori chosen sufficiently large, it can be expected that the iterates $B^\ell$ do never actually reach the box constraints. In this case, no projection is required (cf.~\cref{OPTSYS}). Indeed, in our numerical experiment presented in Section~\ref{sec:NumericalExperiments}, we observe that the partial derivatives of the iterates $B^\ell$ remain bounded uniformly in $\ell$. 
    Therefore, we can argue that suitable box constraints were initially imposed but never reached by the controls during the optimization process. 
    
    For the actual implementation of an $H^1$-projection in a similar situation,
    we refer to \cite[Section 4]{Bartsch2021Numerical}. There, to compute this projection, yet another optimization problem needs to be solved.
\end{rem}

\section{Numerical experiments}
\label{sec:NumericalExperiments}
In this section, we test our implementation of the numerical methods presented in Section~\ref{sec:MC_framework}. 
First, we investigate test cases for our solver of the Vlasov--Poisson system and afterwards, we solve an optimal control problem, where a plasma is to be confined in a certain region in phase space. Our scaling parameters are set to $\mu_x=\mu_v=0.01$ for all of our experiments.

\paragraph*{Landau Damping}
We start with testing the solver for the Vlasov--Poisson model given in \autoref{algo:forward_solver}.
We apply no control (i.e., $B\equiv 0$) and expect to observe the so-called \textit{Landau Damping}, which is a typical feature of collisionless plasma (see, e.g., \cite{Miyamoto2016PlasmaControlledFusion, Chen1984introductionPlasmaFusion, KirchhartWilhelm2023numericalFlowIterationVP}).
Similar to \cite{KirchhartWilhelm2023numericalFlowIterationVP}, we start by considering the following initial distribution of the electrons:
\begin{align}
	\mathring{f}^-(x,v) = \frac{1}{2\pi} \exp\left(-\frac{1}{2}|v|^2\right)(1+ \alpha \cos(k\,x))
	\quad\text{for all }
	(x,v) \in \Omega_x \times \Omega_v.
	\label{eq:ic_Landau}
\end{align}
In this case, we set $x_{\max}=4\pi$, $v_{\max} = 10$, $k=0.5$ and $\alpha = 0.5$. Furthermore, we use $N_f=5\cdot10^4, \Delta t=0.05, \Delta x = 0.1256, \Delta v = 0.1$. Notice that \eqref{eq:ic_Landau} is a perturbation of the equilibrium function given by a Maxwell distribution in velocity and uniform distribution in space.
This initial condition is realized using \autoref{algo:initialize_arbitrary}.
To observe Landau damping, the ions initially have to be in thermal and spatial equilibrium. 
This condition is realized using direct sampling from a uniform distribution in position and a Maxwell distribution in velocity.

As usual, we prescribe periodic boundary conditions for the distribution functions with respect to $x$ at the boundary of $\Omega_x$. These boundary conditions are actually the appropriate choice in accordance with the theoretical results on Landau Damping (see, e.g., the seminal work \cite{Mouhot-Villani}) where the Vlasov--Poisson system is considered in a spatially periodic setup, even though this situation is not covered by the mathematical analysis in the present paper.
Therefore, in this experiment, the exchange of numerical particles over the periodic boundary of $\Omega_x$ is actually intended. 

\begin{figure}
	\centering
	\begin{subfigure}[l]{0.33\textwidth}
		\includegraphics[width=\textwidth]{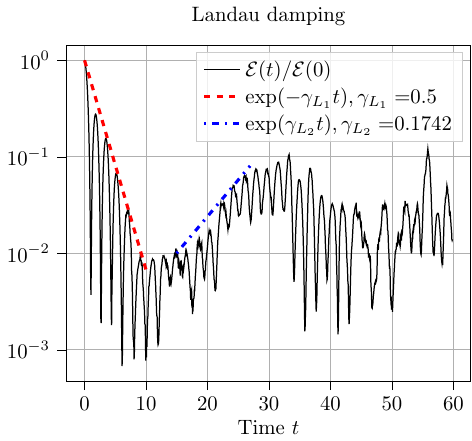}
		\caption{}
		\label{fig:electric_energy}
	\end{subfigure}
	\hspace{-0.01\textwidth}
	\begin{subfigure}[l]{0.33\textwidth}
        \vspace{1cm}
		\includegraphics[width=\textwidth]{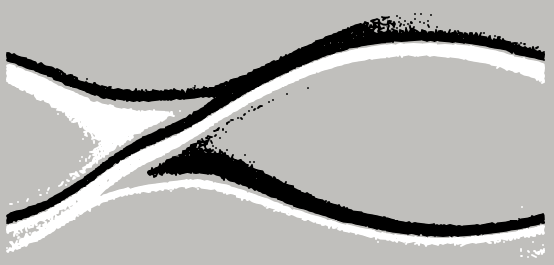}
        \vspace{0.5cm}
		\caption{}
		\label{fig:manifestation_twoStream}
	\end{subfigure}
	\hspace{-0.01\textwidth}
	\begin{subfigure}[l]{0.33\textwidth}
		\includegraphics[width=\textwidth]{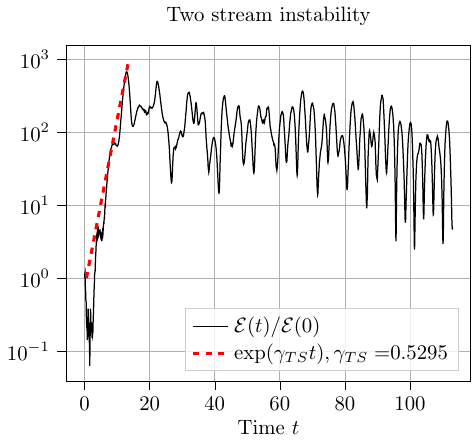}
		\caption{}
		\label{fig:electric_twoStream}
	\end{subfigure}
	\caption{Results of the experiments on Landau damping and two-stream instability.
    (a) Plot of $\mathcal{E}(t)/\mathcal{E}(0)$ (cf. \eqref{eq:electric_energy}) over time with logarithmic scale on ordinate and $\gamma_{L_1}=0.5$ and $\gamma_{L_2}=0.1742$; 
    (b) Manifestation of the two-stream instability. 
    The particles with a positive initial velocity are displayed in black, whereas the particles with a negative initial velocity are displayed in white.
    (c) Evolution of $\mathcal{E}(t)/\mathcal{E}(0)$ over time for the two-stream instability test case with $\gamma_{TS}=0.5295$}.
	\label{fig:TwoStream}	
\end{figure}

In \autoref{fig:electric_energy}, the evolution of the squared norm of the energy is plotted over time.
This means that the evolution of the quantity
\begin{align}
    \mathcal{E}(t) = \int_{\Omega_x} |E_f(t,x)|^2\,\mathrm d x 
    \label{eq:electric_energy}
\end{align}
is displayed.
We expect exponential decay with some rate $\gamma_{L_1}$ which can indeed be observed in \autoref{fig:electric_energy}. 
In particular, our numerical results compare very well to those in \cite[Figure 2]{HoDatta2018HighFidelitySim}. Notice that after a certain time, the value of $\mathcal{E}$ slightly increases again. This is due to the fact that we use a high perturbation with respect to the equilibrium in our initial condition \eqref{eq:ic_Landau}. 
This growth is expected from the nonlinear Landau theory (see, e.g., \cite{BaileyDenavit1970NonlinOscPlasma,Gary19673rdOrderOscillationsPlasma}), 
and even a certain rate $\gamma_{L_2}$ can be calculated (see, e.g., \cite{HoDatta2018HighFidelitySim}).
There, the theoretical values $\gamma_{L_1}=0.562$ and $\gamma_{L_2}=0.168$ are found.
However, as we also encounter numerical noise in our calculation, our growth rates slightly deviate from the theoretical ones. This is due to numerical errors since we only use a first-order particle method (cf.~\cite{EdwardsBridson2012HigOrderPIC}).
We point out that this phenomenon can also be observed in \cite{KirchhartWilhelm2023numericalFlowIterationVP} and  \cite{FinnKnepley2023numericalStudy_LandauDamping_PIC}, and it seems that up to now, such errors cannot be completely avoided.

\paragraph*{Two-stream instability}
The next test case for the forward solver is the appearance of the so-called two-stream instability which is another well-known behavior of plasma (see, e.g., \cite[Chapter 6.6]{Chen1984introductionPlasmaFusion} or \cite{Ghorbanalilu2014PIC_TwoStreamInstab}). We execute this test case in dimension $d=(1,1)$ without any control (i.e., $B\equiv 0$).
The instability occurs if two beams of electrons with opposite directions of velocity evolve over time subject to a static background of ions. 
To this end, the initial distribution of the electrons with respect to velocity is the sum of two Gaussians centered at $v_{\mathrm{TS}}$ and $-v_{\mathrm{TS}}$.
More precisely, we start with an initial condition for the electrons that is slightly perturbed from equilibrium in space: for $(x,v) \in \Omega_x \times \Omega_v$, we choose
\begin{align*}
    \mathring{f}^-(x,v) = 
    \frac{1}{2 \sqrt{2\pi\sigma^2_{\mathrm{TS}}}}\left(
    \exp\left(-\frac{(v-v_{\mathrm{TS}})^2}{2\sigma_{\mathrm{TS}}} \right)
    +
    \exp\left(-\frac{(v+v_{\mathrm{TS}})^2}{2\sigma_{\mathrm{TS}}} \right)
    \right)\left(1+ \alpha \cos(k\,x)\right),
\end{align*}
where $\sigma_{\mathrm{TS}}$ has to be smaller than $v_{\mathrm{TS}}$. As in \cite{HoDatta2018HighFidelitySim}, we consider the values $\sigma_{\mathrm{TS}}=\pi/2, v_{\mathrm{TS}}=0.1, \alpha = 0.01, k = 0.5$. Furthermore, we assume that the ions are Maxwell distributed in velocity and uniformly distributed in space. The numerical parameters are chosen to be $N_f=5\cdot10^4, \Delta t=0.05, x_{\max}=4\pi, v_{\max} = 10, \Delta x = 0.1256, \Delta v = 0.04$. \autoref{fig:manifestation_twoStream} shows the manifestation of the two-stream instability.
On the ordinate, the velocity is plotted and on the abscissa, the position is displayed.
In \autoref{fig:electric_twoStream}, the relative electric energy $\mathcal{E}(t)/\mathcal{E}(0)$ is plotted over time. We point out that our results compare very well with those in \cite{HoDatta2018HighFidelitySim}. In particular, our numerical growth rate $\gamma_{TS} = 0.5295$ does not deviate much from the theoretical rate that is $0.4952$.

\begin{figure}
	\centering
	\begin{subfigure}[l]{0.96\textwidth}
		\includegraphics[width=0.45\textwidth]{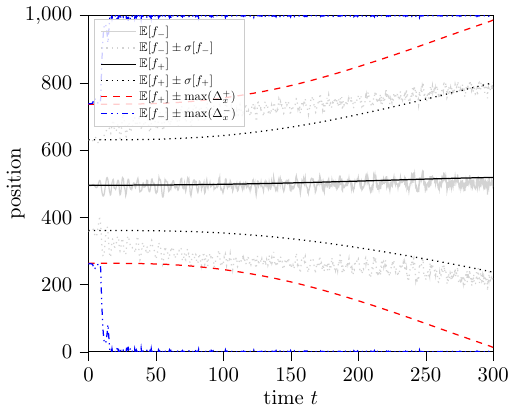}
	\hfill
		\includegraphics[width=0.45\textwidth]{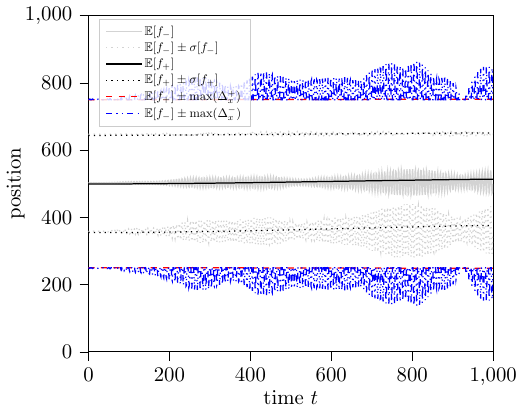}
		\caption{Evolution of mean, variance and maximal deviation from the center of $f^\pm$; see \eqref{eq:max_dev} and \eqref{eq:mean_var} for the definitions. Left: Uncontrolled case; Right: controlled case}
		\label{fig:Evolution_MeanVar_MaxDev_NeuBdry}
	\end{subfigure}
\\
	\begin{subfigure}[l]{0.45\textwidth}
		\includegraphics[width=\textwidth]{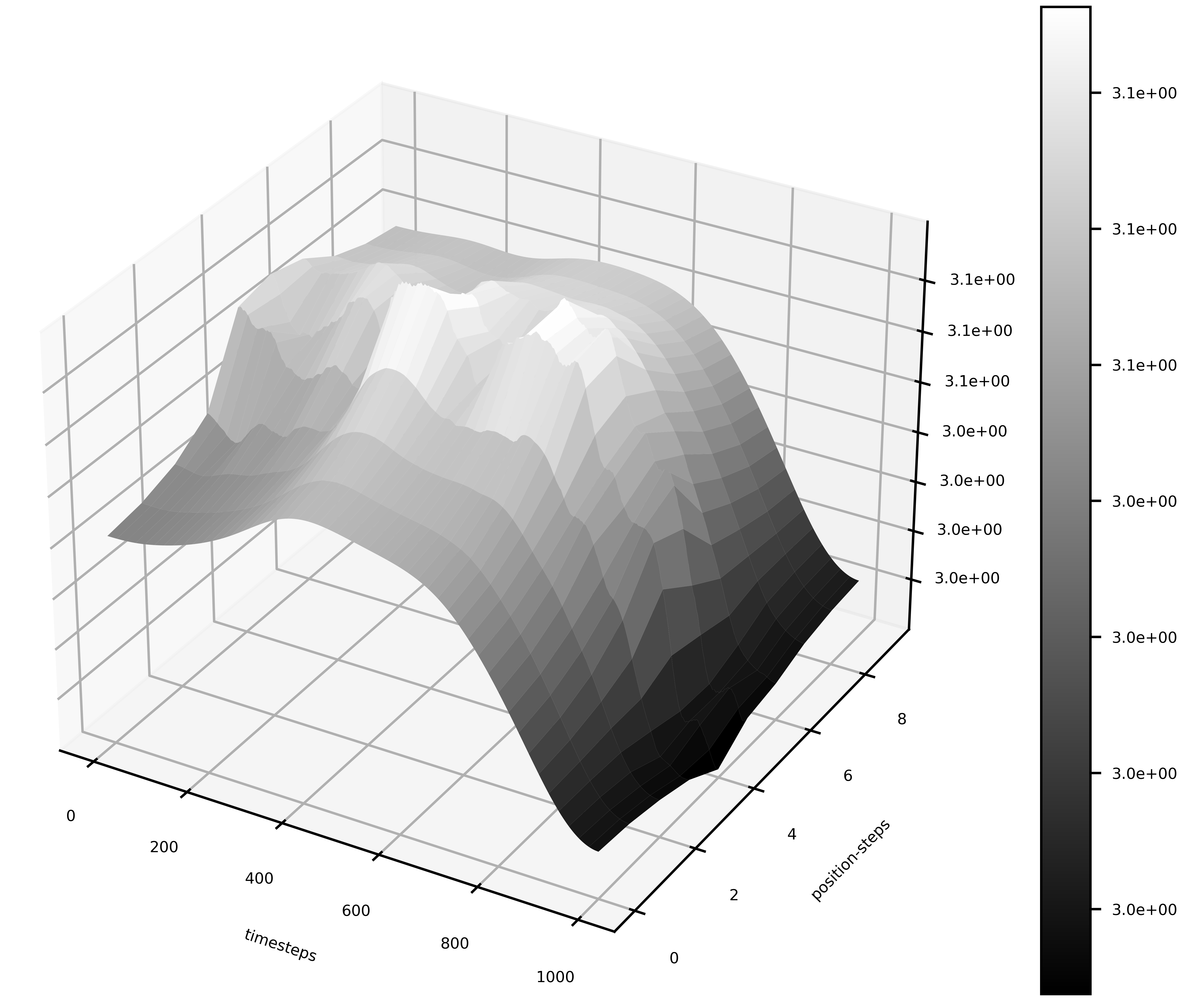}
		\caption{Surface plot of the calculated control}
		\label{fig:Calculated_Control_NeuBdry}
	\end{subfigure}
\hfill
	\begin{subfigure}[l]{0.46\textwidth}
		\includegraphics[width=\textwidth]{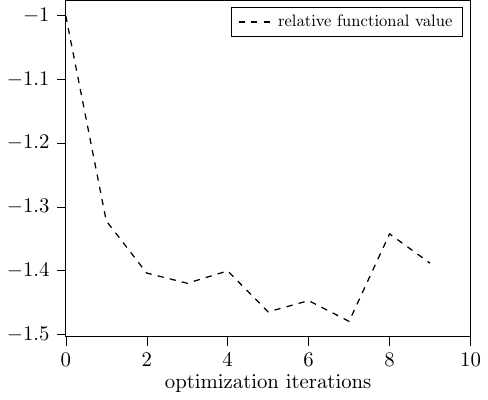}
		\caption{Convergence history of the cost functional $\widetilde J$ (relative values with respect to the value $\widetilde J(B\equiv 0)$)}
		\label{fig:convergence_functional}
	\end{subfigure}
	\caption{Results of optimization procedure given in \cref{algo:OptAlgo} using homogeneous Neumann boundary in time and space for the control}
	\label{fig:Confinement_NeuBdry}
\end{figure}

\paragraph*{Using optimal control to keep a confined initial distribution confined over time}
After these successful tests of the forward solver, we continue with testing our whole optimization method of \autoref{algo:OptAlgo}.
Our goal is to find a control that leads to confinement of the plasma (i.e., the plasma stays away from the boundary of $\Omega_x$) during the time interval $[0,T]$, provided that we start with a confined initial state. More precisely, we use the structure of $\theta^\pm$ and $\varphi^\pm$ given in \eqref{def:phitheta} with the phase-space profiles 
\begin{align*}
	z_{\mathfrak{d}}^\pm = z_T^\pm = \left(x_{\max}/2, 0,0 \right) \, \in \, \RR^{d_x}\times\RR^{d_v},
\end{align*}
which are, in this case, time-independent.
These choices are supposed to ensure that the mean position of the particles is close to the center of the space domain $\Omega_x = (0, x_{\max})$, and that the mean velocity of the particles is close to zero.
As variance matrices we choose
\begin{align*}
	\Sigma_{\theta^\pm} = \Sigma_{\varphi^\pm} =
	\text{diag}(100,\sqrt{\pi/8} \bar{v},\sqrt{\pi/8} \bar{v})
    \, \in \,  \RR^{(d_x+d_v) \times (d_x+d_v)},
\end{align*}
where $\bar{v}=1$ is our chosen expected velocity (cf.~\cite{Garcia00}). Furthermore, we choose our initial data as smooth, bell-shaped functions $\mathring f^\pm$ which are compactly supported in $\Omega_x$ and symmetric about the axis $x=x_{\max}/2$. Specifically, defining $\mathring{\epsilon} = \frac{x_{\max}}{4}$ and $\sigma_v = 0.01$, we choose (cf. \cite[Appendix C.4]{Evans2010PDEs})
\begin{align*}
    \mathring{f}^\pm(x,v) = 
    \begin{cases}
        \dfrac{1}{2\pi \sigma_v^2} \exp\left(
            -\dfrac{2\sigma_v^2}{|v|^2}
        \right)
        \exp \left(
            -\dfrac{1}{1- (x/\mathring{\epsilon})}
        \right)
        & \text{if } \left| x - \frac{x_{\max}}{2} \right| \leq \mathring{\epsilon},
        \\
        0 & \text{else}.
    \end{cases}    
\end{align*}
The numerical parameters are set to be $N_f=10^3$, $\Delta t = 0.3$, $x_{\max} = 10^3$, $\Delta x = 10^2$, $v_{\max}=5\cdot10^3$, $\Delta v = 10^3$, and the control weights are set to $\alpha = 10^{-4}$ and $\kappa_t = \kappa_x = 10^{-4}$. In the left figure in \autoref{fig:Evolution_MeanVar_MaxDev_NeuBdry}, the results of the evolution of the model without applying any control are plotted.
To be precise, the evolution of the mean and variance of the distribution of electrons and ions are displayed over $[0,T]$ together with the maximal deviation $\max(\Delta_x^\pm)$ of the numerical particles from the center of the spatial domain $\Omega_x$.
For every timestep $k = 1,\ldots,N_t$, this maximal deviation is defined as
\begin{align}
	\max(\Delta_x^\pm) \coloneqq \max_{p \in {1,\ldots,N_f}} \left( \left|F_{\pm}^k[p].x - \frac{x_{\max}}{2}\right| \right).
	\label{eq:max_dev}
\end{align} 
The (discrete) mean and the variance with respect to position are defined as usual, that is,
\begin{align}
	\EE[f](t^k) \coloneqq \frac{1}{N_f}\sum_{p=1}^{N_f} F^k[p].x,
	\qquad
	\mathbb{\sigma}^2[f](t,x) \coloneqq \frac{1}{N_f} \left(
	\sum_{p=1}^{N_f} (F^k[p].x)^2 - (\EE[f](t^k))^2\right).
	\label{eq:mean_var}
\end{align}
It becomes clear from the left figure in \autoref{fig:Evolution_MeanVar_MaxDev_NeuBdry} that the support of $f^\pm$ is growing and eventually fills the complete domain $\Omega_x$ if no control is applied.
On the contrary, as shown in \autoref{fig:Evolution_MeanVar_MaxDev_NeuBdry} on the right, the locally optimal control $B$ detected by our optimization method leads to the effect that the supports of $f^\pm$ stay away from the boundary of the interval $\Omega_x$. 
In fact, the support of the distribution function $f^+$ (that is associated with the ions) hardly changes over the course of time. 
The calculated optimal control is depicted in \cref{fig:Calculated_Control_NeuBdry}.
It is worth mentioning that we started our optimization method without assuming further knowledge of the control, i.e., we initialized the control by simply choosing $B\equiv 0$. 
Altogether, this illustrates the ability of our optimal control approach to find an external magnetic field $B$ that keeps the confined initial state confined over the course of time, at least on a short time interval.

In \cref{fig:convergence_functional}, we plot the convergence history of the functional. 
Recall that we use a particle method and re-sample the initial condition every time when the state equations have to be solved. In this way, we intend to find a (locally) optimal control that does not only work well with a specific initial distribution of the numerical particles.
However, this then leads to a non-monotonic behaviour of the convergence of the cost functional, which becomes oscillating towards the end of the convergence.

\begin{figure}
	\centering
	\begin{subfigure}[l]{0.98\textwidth}
		\includegraphics[width=.48\textwidth]{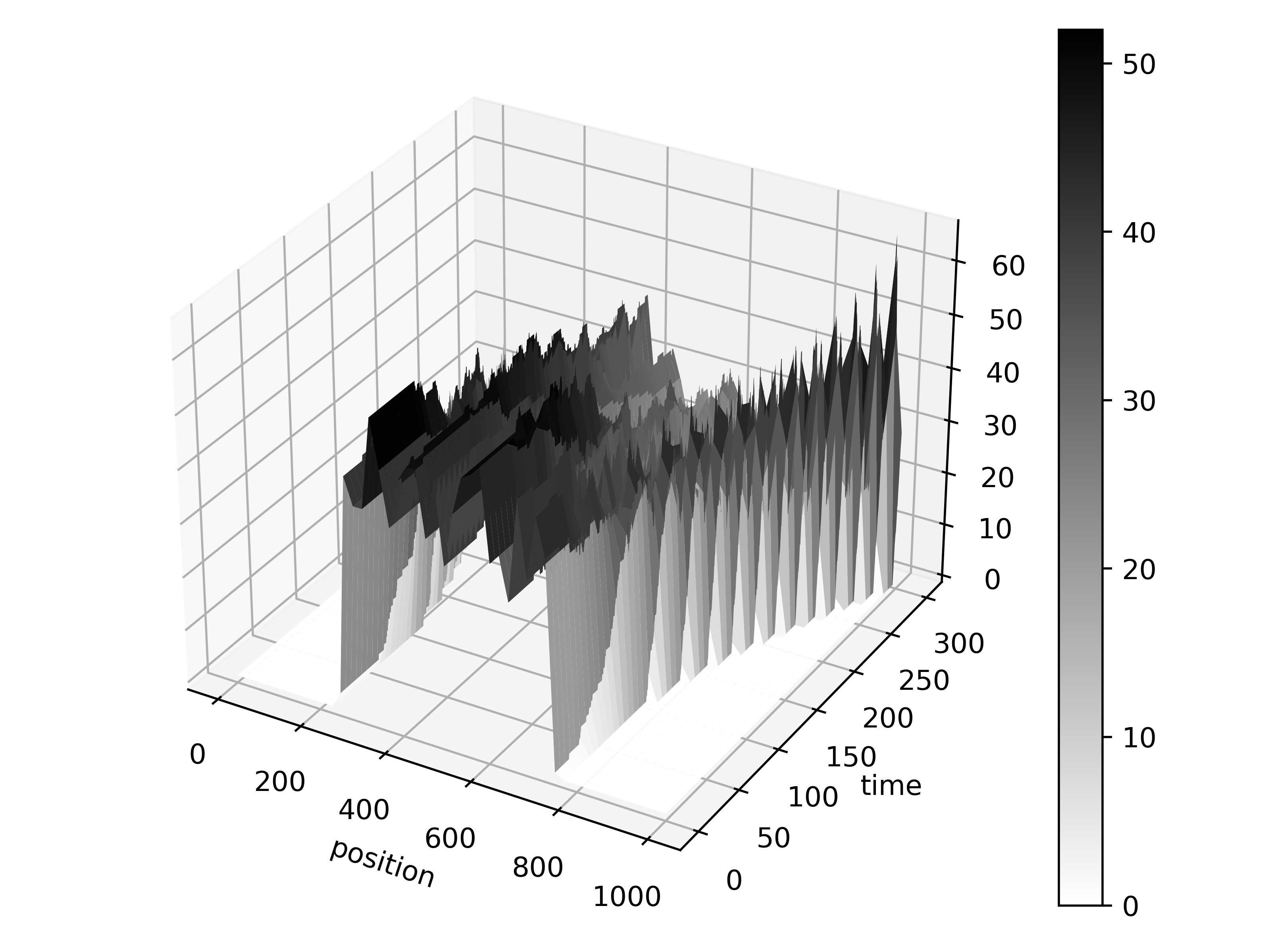}
	\hfill
		\includegraphics[width=.48\textwidth]{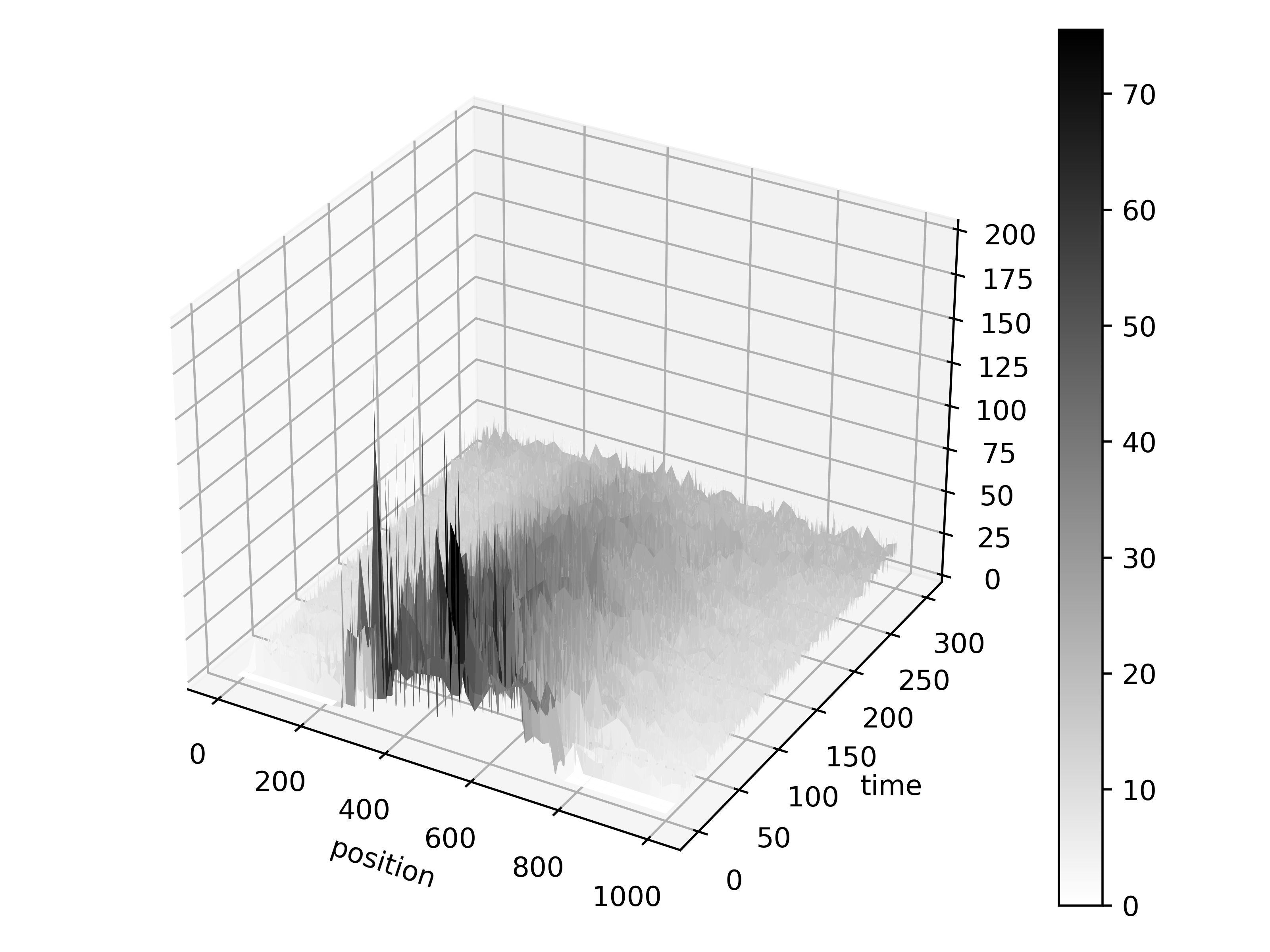}
		\caption{Uncontrolled densities (left: ions, right: electrons)}
		\label{fig:Uncontrolled_densities}
	\end{subfigure}
	\begin{subfigure}[l]{0.98\textwidth}
		\includegraphics[width=.48\textwidth]{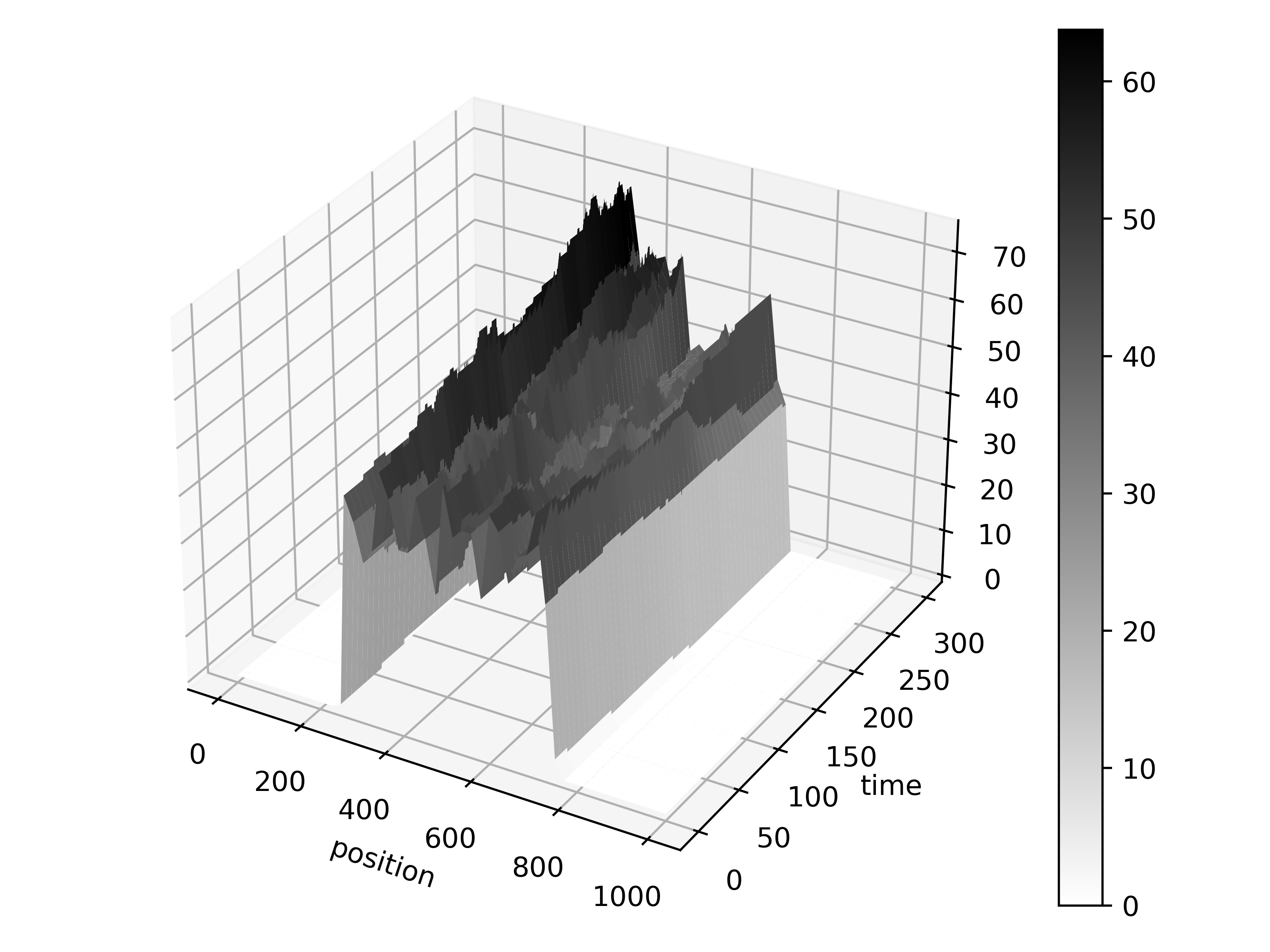}
	\hfill
		\includegraphics[width=.48\textwidth]{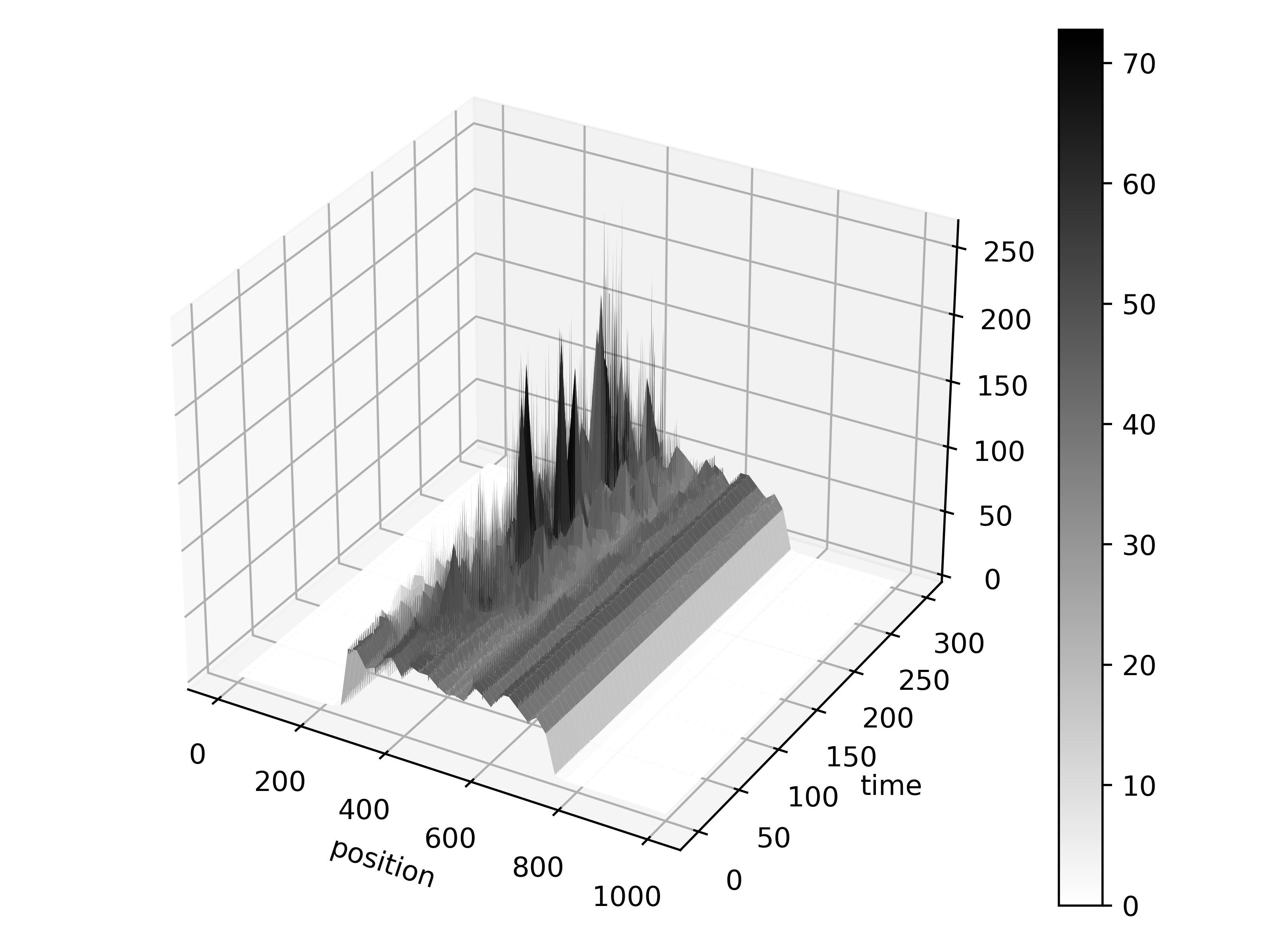}
		\caption{Controlled densities (left: ions, right: electrons)}
		\label{fig:Controlled_densities}
	\end{subfigure}
	\caption{Evolution of numerical ion (left) and electron (right) charge densities (cf.~\eqref{eq:rho_ions_electrons}). 
 Above: behaviour when applying no control; Below: behaviour after applying the control calculated by our optimization procedure that is depicted in \cref{fig:Calculated_Control_NeuBdry}.}
	\label{fig:density_plots}
\end{figure}

In \cref{fig:density_plots}, we plot the densities 
\begin{align}
    \rho_{f^\pm}(t,x) \coloneqq
    \int_{\RR^{d_v}} f^\pm(t,x,v) \, \mathrm{d}v,
    \label{eq:rho_ions_electrons}
\end{align}
using no control ($B\equiv 0$) and the control calculated using our optimization strategy.
These plots also nicely illustrate that our control is able to confine the particles away from the boundary. 
In \cref{fig:Uncontrolled_densities}, which shows the uncontrolled scenario, the electrons reach the boundary after a very short time, whereas the ions reach the boundary towards the end of the time interval.
On the contrary, in \cref{fig:Controlled_densities}, the particles clearly stay away from the boundary during the whole time interval.

 \section{Acknowledgements}
Stefania Scheurer was partially funded by the Zukunftskolleg of the University of Konstanz.
Jan Bartsch was partially supported by the SFB 1432 ``Fluctuations and Nonlinearities'' (Project-ID 425217212) of the Deutsche Forschungsgemeinschaft (DFG, German Research Foundation) and
Patrik Knopf was partially supported by the RTG 2339 ``Interfaces, Complex Structures, and Singular Limits'' of the DFG.
Moreover, Jörg Weber was supported by the Swedish Research Council (grant no 2020-00440).
All of this support is gratefully acknowledged.
We also thank Leo Basov (DLR Göttingen) for his very helpful comments concerning the numerical validation of the implementation.  
We would also like to express our gratitude to the anonymous referees for their helpful questions and remarks.

%%%%% References
%\bibliographystyle{acm}
%\bibliography{literature_PlasmaMonteCarlo}{}
%%%%%

\end{document}